\documentclass[10pt,a4paper,reqno]{amsart}

\usepackage[british]{babel}
\usepackage{csquotes}
\usepackage{epsfig}
\usepackage{amsmath}
\usepackage{bbm}
\numberwithin{equation}{section}
\usepackage[hidelinks]{hyperref}

\usepackage[left=2.5cm,right=2.5cm,top=2.5cm,bottom=2.5cm,footskip=1.5cm,headsep=1cm]{geometry}
\usepackage{appendix}

\newcommand{\interior}[1]{%
  {\kern0pt#1}^{\mathrm{o}}%
}

\usepackage[T1]{fontenc} 
\usepackage{lmodern} 
\rmfamily 
\DeclareFontShape{T1}{lmr}{b}{sc}{<->ssub*cmr/bx/sc}{}
\DeclareFontShape{T1}{lmr}{bx}{sc}{<->ssub*cmr/bx/sc}{}
\usepackage{lipsum}
\usepackage{mathtools}
\usepackage{amsmath}
\usepackage{amssymb}
\usepackage{tikz}
\usetikzlibrary{matrix,positioning,decorations.pathreplacing,calc}
\usetikzlibrary{
decorations.text,%
decorations.markings,%
shadows}
\usepackage{adjustbox}
\usepackage{graphicx} 

\usepackage{caption}
\usepackage{subcaption}

\usepackage{dsfont} 
\usepackage[format=hang]{caption}

\usepackage[normalem]{ulem}

\usepackage{bm}

\usepackage[
sorting=nty,
maxnames=99,
maxalphanames=5,
natbib=true,
sortcites]{biblatex}

\DeclareNameAlias{default}{family-given}

\graphicspath{{figures/}}

\AtEveryBibitem{
 \clearfield{url}
 \clearfield{issn}
 \clearfield{isbn}
 \clearfield{urldate}
 
 \ifentrytype{book}{
  \clearfield{pages}}{%
 }
}

\usepackage{xargs}                      
\usepackage[prependcaption,textsize=tiny,textwidth=3cm]{todonotes}
\newcommandx{\unsure}[2][1=]{\todo[linecolor=red,backgroundcolor=red!25,bordercolor=red,#1]{#2}}
\newcommandx{\change}[2][1=]{\todo[linecolor=blue,backgroundcolor=blue!25,bordercolor=blue,#1]{#2}}
\newcommandx{\info}[2][1=]{\todo[linecolor=OliveGreen,backgroundcolor=OliveGreen!25,bordercolor=OliveGreen,#1]{#2}}
\newcommandx{\improvement}[2][1=]{\todo[linecolor=black,backgroundcolor=black!25,bordercolor=black,#1]{#2}}
\newcommandx{\thiswillnotshow}[2][1=]{\todo[disable,#1]{#2}}
\usepackage{amsthm}
\usepackage[noabbrev,capitalize]{cleveref}
\crefname{proposition}{Proposition}{Propositions}
\crefname{equation}{}{}

\newtheorem{theorem}{Theorem}[section]
\newtheorem{lemma}[theorem]{Lemma}
\newtheorem{proposition}[theorem]{Proposition}
\newtheorem{corollary}[theorem]{Corollary}

\theoremstyle{definition}

\newtheorem{remark}[theorem]{Remark}

\crefname{assumption}{Assumption}{Assumptions}
\crefname{definition}{Definition}{Definitions}
\crefname{corollary}{Corollary}{Corollaries}
\crefname{enumi}{item}{items}

\usepackage{fancyhdr}

\fancypagestyle{plain}{%
\fancyhead[C]{}
\fancyfoot[C]{}

}
\fancypagestyle{nosection}{%
\fancyhead[CE]{}
\fancyhead[CO]{}
\fancyfoot[CE,CO]{\thepage}
}
\newcommand{\nm}{\noalign{\smallskip}}
\newcommand{\ds}{\displaystyle}

\DeclareMathOperator{\R}{\mathbb{R}}

\renewcommand{\tilde}{\widetilde}
\renewcommand{\hat}{\widehat}

\renewcommand{\qedsymbol}{$\blacksquare$}

\usepackage{fancyhdr}

\fancypagestyle{plain}{%
\fancyhead[C]{}
\fancyfoot[C]{}

}
\fancypagestyle{nosection}{%
\fancyhead[CE]{}
\fancyhead[CO]{}
\fancyfoot[CE,CO]{\thepage}
}

\renewcommand{\epsilon}{\varepsilon}

\renewcommand{\tilde}{\widetilde}
\renewcommand{\hat}{\widehat}

\usepackage{tcolorbox}
\usetikzlibrary{tikzmark,calc}

\usepackage{enumerate}
\usepackage{upgreek}

\theoremstyle{remark}

\numberwithin{equation}{section}

\usepackage{setspace}

\def\1{\textbf{\rm 1}}

\makeatletter
\@namedef{subjclassname@2020}{\textup{2020} Mathematics Subject Classification}
\makeatother

\begin{document}

\title{Sharp Gaussian Brunn-Minkowski inequalities}
\author[B. Ammari]{Bader Ammari}
\address[]{Department of Mathematics, ETH Z\"{u}rich, R\"{a}mistrasse 101, CH-8092 Z\"{u}rich, Switzerland}
\email{bader.ammari@math.ethz.ch}
\begin{abstract}
    We prove Gaussian Brunn-Minkowski inequalities for multiple sets. These are of the type: $$\prod_{i=1}^I \gamma_{\Sigma_i}(K_i)^{\tilde \lambda_i} \le C_{\text{sharp}}\gamma_\Gamma\Big(\sum_{i=1}^I \lambda_i K_i\Big),$$ with convex sets $K_i \subset \mathbb{R}^n$, whose barycenter $b_{K_i,\Sigma_i}$ with respect to the Gaussian measure $\gamma_{\Sigma_i}$ is zero, where $\lambda_i, \tilde \lambda_i \in (0;1)$ such that $\sum_{i=1}^I \lambda_i=\sum_{i=1}^I \tilde \lambda_i=1.$ In addition, we provide natural conditions on the matrices $\Gamma, \Sigma_1,\cdots,\Sigma_I \in \mathrm{Sym}_{>0}(n)$, which are given through a pointwise inequality involving the Karcher mean. Furthermore, we study equality cases of the corresponding functional forms in such situations. We also investigate the role of centering in the (B)-property by giving a counterexample for convex bodies with zero Lebesgue barycenters in dimensions at least three and thereby answering negatively a question of Faifman-Polterovich in \cite{Faifman_Polt}. Finally, a new slicing approach is introduced, which allows us to pass from centered convex sets to sets with arbitrary barycenters such that their weighted sum is convex. Based on this idea, we show a new Gaussian Brunn-Minkowski inequality $$e^{|\sum_{i=1}^I \lambda_i b_{K_i,I_n}|^2/2n} \gamma_n \Big(\sum_{i=1}^I \lambda_i K_i \Big)^{1/n} \ge \sum_{i=1}^I \lambda_i e^{|b_{K_i,I_n}|^2/2n} \gamma_n(K_i)^{1/n},$$ and consider a further application in the direction of the Gaussian correlation inequality, together with the rigidity for both inequalities. 
\end{abstract}
\maketitle

\vspace{0.2cm}

\vspace{0.2cm}
\vspace{0.5cm}

\section{Introduction}

The celebrated Brunn-Minkowski inequality states that for all $K,L \subset \mathbb{R}^n$ measurable sets with positive measure, $$|K+L|^{1/n} \ge |K|^{1/n}+|L|^{1/n}.$$ By homogeneity of the Lebesgue measure, it can be written as $$|(1-t)K+tL| \ge |K|^{1-t}|L|^t \quad \forall t \in [0;1].$$ Its well-known functional form is the Pr\'ekopa-Leindler inequality, which says that $$\int_{\mathbb{R}^n}\sup_{z=t x+(1-t)y}f(x)^t g(y)^{1-t}\, dz \ge \Big(\int_{\mathbb{R}^n} f(x)\, dx\Big)^t \Big(\int_{\mathbb{R}^n}g(y)\, dy \Big)^{1-t},$$ for any $f,g: \mathbb{R}^n \rightarrow [0;+\infty)$ with finite and positive integrals. The previous functional inequality is known thanks to the pioneering work of Barthe \cite{Barthe} to be a reverse Brascamp-Lieb inequality. We refer the reader to the beautiful survey of \cite{Gardner_BM} and the references therein.

\hfill \\
Brascamp-Lieb inequalities were first studied in \cite{Brascamp_Lieb} and are functional inequalities whose sharp constant is completely determined by centered Gaussian functions (commonly known as the Gaussian saturation phenomenon) as proven by Lieb in \cite{Lieb}. These inequalities encompass many important inequalities such as H\"older, Young convolution or Loomis-Whitney and immediately found important applications in the field of convex geometry thanks to the works of Keith Ball (see \cite{KBall,KBall2}). 

\hfill \\
Two natural and closely related questions with regard to such functional or geometric inequalities can then be addressed: 
\begin{itemize}
\item[(i)] By imposing natural assumptions on the sets or functions considered, could certain Brunn-Minkowski type inequalities become true? 
\item[(ii)] Can we replace the Lebesgue measure by other kinds of measures for which the statement would still hold, up to the aforementioned assumptions?
\end{itemize}
These two questions turn out to be highly general, very challenging and have motivated a major direction of research in the field of convex geometry. A first natural assumption on sets (resp. functions) is convexity (resp. log-concavity) but symmetry with respect to the origin ($K=-K$ for a set and requiring evenness for functions) proves to be also a crucial assumption.

\hfill \\
The usual measures considered for the validity of such inequalities are log-concave measures and more specifically, Gaussian measures. Recall that for $\Sigma \in \mathrm{Sym}_{>0}(n)$, we write $\gamma_\Sigma$ for the centered Gaussian probability measure with covariance $\Sigma$ and $$\gamma_\Sigma:=\frac1{\sqrt{\det(2\pi \Sigma)}}e^{-\frac12 \langle \Sigma^{-1}x,x\rangle}\, dx.$$ 

\hfill \\
A first example of such inequality is linked to the (B)-theorem, a problem proposed by Banaszczyk and popularized by \L atala in \cite{Latala}, then subsequently solved by Cordero-Erausquin, Fradelizi and Maurey in \cite{Cordero_Fradelizi_Maurey}, where it is proven that if $K \subset \mathbb{R}^n$ is a symmetric convex set, then the function $$(t_1,\cdots,t_n) \mapsto \gamma_n(e^{\mathrm{diag}(t_1,\cdots,t_n)}K) \text{ is log-concave on } \mathbb{R}^n.$$ The proof relies on Poincar\'e-type inequalities for $\gamma_K$, the normalized restriction of $\gamma_n$ to the symmetric convex set $K$, and the so-called H\"ormander's $L^2$ method. 

\hfill \\
This method has become crucial to study Brunn-Minkowski type inequalities for the Gaussian measure and to understand (B)-properties of log-concave measures. First, by the Pr\'ekopa-Leindler inequality, for all $K,L \subset \mathbb{R}^n$ with positive Gaussian measure, $$\gamma_n((1-t)K+tL) \ge \gamma_n(K)^{1-t}\gamma_n(L)^t \quad \forall t \in [0;1],$$ however, as homogeneity fails for the Gaussian measure, this does not recover an additive-type inequality of the form \begin{equation} \label{Basic_Gaussian_BM} \gamma_n((1-t)K+tL)^{1/n} \ge (1-t)\gamma_n(K)^{1/n}+t\gamma_n(L)^{1/n} \quad \forall t \in [0;1].\end{equation} In fact, such a statement was conjectured by Gardner and Zvavitch in \cite{Gardner_Zvavitch} for convex sets $K,L$ containing the origin but a counterexample was provided in \cite{Nayar_Tkocz}. Thanks to the works of Kolesnikov-Livshyts \cite{Kolesnikov_Livshyts} and later of Eskenazis-Moschidis \cite{Eskenazis_Moschidis}, \cref{Basic_Gaussian_BM} was proven for symmetric convex sets and later, Cordero-Erausquin and Rotem in \cite{Cordero_Rotem} showed that such a Brunn-Minkowski inequality holds for all rotationally invariant log-concave measures. Additionally, they proved that the (B)-property, \textit{i.e.} for every symmetric convex $K \subset \mathbb{R}^n$ and every symmetric matrix $A$, $$t \mapsto \mu(e^{tA}K) \text{ is log-concave }$$ holds for any rotationally invariant log-concave measure $\mu.$  

\hfill \\
For more recent works using such approaches, we refer the reader to \cite{Cordero_Eskenazis,Cordero_Rotem_Livshyts}. Let us note that in the case of the standard Gaussian measure, a proof of the (B)-property (or in other words, the strong form of the Gaussian (B) inequality) was already fully shown by Brascamp and Lieb in \cite{Brascamp_LiebBThm} and employs the rotation trick, a crucial argument in the proof of the Gaussian saturation, together with the variance Brascamp-Lieb inequality applied to a linear odd function (see the appendix of \cite{Cordero_Eskenazis}). 

\hfill \\
In fact, only focusing on the (B)-theorem for Gaussian measures with \textit{arbitrary} covariances is linked to the fundamental log Brunn-Minkowski problem, first formulated by Böröczky, Lutwak, Yang and Zhang in \cite{Borockzy_Lutwak_Yang_Zhang}, connected to important problems of uniqueness of the cone measures and Monge-Amp\`ere, see \cite{Borockzy_Lutwak_Yang_Zhang2} and recently announced to be solved by OpenAI in \cite{OpenAI_logBM}. By \cite{Saroglou}, the log Brunn-Minkowski inequality, which is yet an example of pivotal improvement compared to the Brunn-Minkowski inequality under convexity and evenness, $$\Big|\bigcap_{u \in \mathbb{S}^{n-1}} \{x \in \mathbb{R}^n: \langle x,u\rangle \le h_K(u)^{1-t}h_L(u)^t \} \Big| \ge |K|^{1-t}|L|^t,$$ where $h_K(x):=\sup_{y \in K}\langle x,y\rangle$, is equivalent to the fact that for all $\Sigma \in \mathrm{Sym}_{>0}(n)$ and for every diagonal matrix $A \in \mathrm{Sym}(n)$, $$t \mapsto \gamma_\Sigma(e^{tA}[-1;1]^n) \text{ is log-concave.}$$
However, in \cite{Cordero_Rotem2}, it is shown that $t \mapsto \gamma_\Sigma(e^{tA}K)$ may not be necessarily log-concave for a symmetric convex body $K \subset \mathbb{R}^n$ and a diagonal matrix $A.$ 

\hfill \\
Lastly, another recent way of proving \cref{Basic_Gaussian_BM} has appeared in \cite{Aishwarya_Li} using information-theoretic tools and constructing a transport coupling, which satisfies the entropic inequality corresponding to the Gaussian Brunn-Minkowski inequality. Building on \cite{Aishwarya_Rotem}, Han-Zhu in \cite{Han_Zhu} proved that the displacement interpolation between two $1$-strongly log-concave probability measures has also Poincar\'e constant at most $1$, in particular settling \cite[Question 1.6]{Aishwarya_Rotem}, and deduced a Gaussian Brunn-Minkowski inequality for two sets: for $K_1,K_2 \subset \mathbb{R}^n$ convex sets with $\gamma_n(K_1),\gamma_n(K_2)>0$, \begin{equation} \label{Han_Zhu_BM} e^{|(1-t)b_{K_1}+tb_{K_2}|^2/2n}\cdot \gamma_n\Big((1-t)K_1+tK_2\Big)^{1/n} \ge (1-t)e^{|b_{K_1}|^2/2n}\gamma_n(K_1)^{1/n}+te^{|b_{K_2}|^2/2n}\gamma_n(K_2)^{1/n},\end{equation} where $$b_{K_i}:=\frac1{\gamma_n(K_i)}\int_{K_i} x\, d\gamma_n(x) \text{ for } i=1,2.$$ 

\hfill \\
Let us now come back to the relationship with Brascamp-Lieb inequalities. As explained previously, these functional inequalities have important consequences but as they hold for any functions, it is again worth considering whether by restricting functions to for example log-concave and even functions, one could have interesting inequalities. Nakamura and Tsuji in \cite{Shohei_Hiroshi} introduced such a viewpoint of improvement under evenness and log-concavity in the theory of Brascamp-Lieb inequalities, first through inverse Brascamp-Lieb inequalities considered in the general case by \cite{Barthe_Wolff}. Later, Milman in \cite{Milman} realized that such a formulation yielded the celebrated Gaussian correlation inequality (see \cite{Royen,Latala_Matlak} for the original arguments) and Nakamura and Tsuji in \cite{Shohei_Hiroshi2} pushed further this idea to obtain a generalized Gaussian correlation inequality with multiple centered convex sets.
These works culminated in \cite{Strong_Gaussian_Correlation} to obtain a general Gaussian Forward-Reverse Brascamp-Lieb (FRBL) with the required assumptions of log-concavity and centering (the Lebesgue case was known from \cite{Courtade_Liu}), which allowed to prove the strong Gaussian correlation inequality. 

\hfill \\
Our main contributions in this paper are as follows: we first determine the sharp constants in Gaussian Brunn–Minkowski inequalities for arbitrary finite families of convex sets centered with respect to possibly different Gaussian measures. We characterize equality under Karcher mean conditions and in the associated functional interpolation inequalities. We also introduce a slicing method that yields Gaussian inequalities with barycenter corrections, including Brunn-Minkowski inequalities for measurable sets whose weighted Minkowski sum is convex and correlation inequalities for convex sets with a common Gaussian barycenter. Finally, we investigate the role of centering in the (B)-property, constructing counterexamples and showing that optimization over translations restores log-concavity.

\hfill \\
We would like to note that \cref{Sharp_Constant_BM_Ineq} was obtained independently of \cite{Heterovariant_BM}. In the two-input case, \cref{Sharp_Constant_BM_Ineq} coincides with the sharp heterovariant Gaussian Brunn-Minkowski inequality established by Milman, Nakamura and Tsuji in the second version of \cite{Heterovariant_BM}. Both works employ the Gaussian Forward-Reverse Brascamp-Lieb theorem of Milman, Nakamura and Tsuji \cite{Strong_Gaussian_Correlation} to reduce the sharp-constant problem to a matrix optimization one. Here, the optimization is subsequently treated for arbitrary finite families using Schur complements and induction on the dimension. The results below also provide functional inequalities and equality characterizations under the Karcher mean conditions.

\subsection{Gaussian inequalities, convexity and centering}
\hfill \\
Let $I \ge 2$ and $\lambda=(\lambda_1,\cdots,\lambda_I) \in \Delta_I^\circ$, where the relative interior of the unit simplex is denoted by $$\Delta_I^\circ=\Big\{(\lambda_1,\cdots,\lambda_I): \lambda_i>0, \quad \sum_{i=1}^I \lambda_i=1 \Big\}.$$

\hfill \\
If $a_i>0$ for $1 \le i \le I$, denote their weighted power mean by $$\mathcal M_p^\lambda(a_1,\cdots,a_I)=\begin{cases}
    \ds \Big(\sum_{i=1}^I \lambda_i a_i^p\Big)^{1/p}, \quad 0<p<+\infty, \\
    \ds \prod_{i=1}^I a_i^{\lambda_i}, \quad p=0, \\
    \ds \max_{1 \le i \le I}a_i, \quad p=+\infty,
\end{cases}$$ for any $p \in [0;+\infty]$ and $\lambda \in \Delta_I^\circ.$ Moreover, the mean is set to zero if any $a_i$ is zero. 

\hfill \\
We will write for a given matrix $\Gamma \in \mathrm{Sym}(\mathbb{R}^n)$, $$\Gamma_E=P_E\Gamma|_E,$$ where $P_E$ denotes the orthogonal projection onto a subspace $E \subset \mathbb{R}^n.$

\hfill \\
Our first main result is as follows.

\begin{theorem} \label{Sharp_Constant_BM_Ineq}
Let $\lambda, \tilde \lambda \in \Delta_I^\circ$ and let $\Gamma, \Sigma_i \in \mathrm{Sym}_{>0}(n).$ Then, for any convex sets $K_1,\cdots,K_I \subset \mathbb{R}^n$ such that $$\int_{K_i}x\, d\gamma_{\Sigma_i}(x)=0 \quad \forall 1 \le i \le I,$$ one has $$\prod_{i=1}^I \gamma_{\Sigma_i}(K_i)^{\tilde \lambda_i} \le \max_{E \text{ linear subspace of } \mathbb{R}^n} \prod_{i=1}^I \Big(\frac{\tilde \lambda_i}{ \lambda_i}\Big)^{\dim(E) \tilde \lambda_i}\Bigg[ \frac{\det(\Gamma_E)}{\prod_{i=1}^I (\det((\Sigma_i)_E))^{\tilde \lambda_i}}\Bigg]^{1/2} \gamma_\Gamma \Big(\sum_{i=1}^I \lambda_i K_i\Big).$$ 
\end{theorem}

\hfill \\

\begin{remark}
    Note that the constant is sharp. Indeed, by fixing a subspace $E$ of dimension $\dim(E) \ge 1$, we choose a symmetric convex body $B \subset E$ with positive $\dim(E)$-Lebesgue measure and set $$K_i^{(\epsilon)}=\frac{\tilde \lambda_i}{\lambda_i}\epsilon B+E^\perp.$$ On the one hand, $$\sum_{i=1}^I \lambda_i K_i^{(\epsilon)}=\epsilon B+E^\perp,$$ so that $$\gamma_\Gamma\Big(\sum_{i=1}^I \lambda_i K_i^{(\epsilon)}\Big) \sim \epsilon^{\dim(E)}\frac{|B|}{(2\pi)^{\dim(E)/2}\sqrt{\det \Gamma_E}} \text{ as } \epsilon \rightarrow 0.$$  On the other hand, $$\gamma_{\Sigma_i}(K_i^{(\epsilon)})=\gamma_{(\Sigma_i)_E}\Big(\frac{\tilde \lambda_i}{\lambda_i}\epsilon B\Big) \sim \Big(\frac{\tilde \lambda_i}{\lambda_i}\Big)^{\dim(E)} \epsilon^{\dim(E)}\frac{|B|}{(2\pi)^{\dim(E)/2}\sqrt{\det (\Sigma_i)_E}} \text{ as } \epsilon \rightarrow 0.$$ Hence, $$\frac{\prod_{i=1}^I \gamma_{\Sigma_i}(K_i^{(\epsilon)})^{\tilde \lambda_i}}{\gamma_\Gamma\Big(\sum_{i=1}^I \lambda_i K_i^{(\epsilon)}\Big)} \rightarrow_{\epsilon \downarrow 0} \prod_{i=1}^I \Big(\frac{\tilde \lambda_i}{ \lambda_i}\Big)^{\dim(E) \tilde \lambda_i}\Bigg[ \frac{\det(\Gamma_E)}{\prod_{i=1}^I (\det((\Sigma_i)_E))^{\tilde \lambda_i}}\Bigg]^{1/2}.$$ The term corresponding to $E=\{0\}$ is $1$ and attained by taking all $K_i=\mathbb{R}^n.$
\end{remark}

\hfill \\
Furthermore, there are substantial cases of interest where the sharp constant has an explicit form, which leads to interpolated forms of inequalities between the heterovariant Gaussian (B)-theorem and the heterovariant Gaussian Brunn-Minkowski inequality. 

\hfill \\
Define the Karcher mean matrix $M$ of $R_1,\cdots,R_I \in \mathrm{Sym}_{>0}(n)$, denoted by $\mathcal O_\lambda(R_1,\cdots,R_I)$, as the unique positive definite solution of $$\sum_{i=1}^I \lambda_i \log(M^{-1/2}R_iM^{-1/2})=0.$$ 
\hfill \\ For special properties of the Karcher mean, we refer the reader to \cite{Lawson_Lim}. 
In the case of two matrices, the Karcher mean becomes the familiar geometric mean, \textit{i.e.} recall that for $A,B \in \mathrm{Sym}_{>0}(n)$ and $t \in [0;1]$, the $t$-geometric mean of $A$ and $B$ is defined by $$A \#_t B=A^{1/2}\big(A^{-1/2}BA^{-1/2}\big)^t A^{1/2}.$$ Indeed, for $t \in (0,1),$ if $$(1-t)\log(M^{-1/2}AM^{-1/2})+t \log(M^{-1/2}BM^{-1/2})=0,$$ we may set $C=M^{-1/2}AM^{-1/2}$ and $D=M^{-1/2}BM^{-1/2}$, so that $$\log D=-\frac{1-t}t\log C \Longrightarrow D=C^{-(1-t)/t} \Longrightarrow C\#_t D=C^{1-t}D^t=I_n.$$ Then, $$A\#_t B=(M^{1/2}CM^{1/2})\#_t (M^{1/2}DM^{1/2})=M^{1/2}(C\#_t D)M^{1/2}=M.$$  

\hfill \\

\begin{theorem} \label{Main_Thm_IneqPart}
    Let $p \in [0;1], \lambda, \tilde \lambda \in \Delta_I^\circ, S_i \in \mathrm{Sym}(n)$ for any $1 \le i \le I$ and $\Gamma, \Sigma_i \in \mathrm{Sym}_{>0}(n).$ 
    If $p=0$, assume that $$\Gamma \le \mathcal O_\lambda(e^{S_1}\Sigma_1e^{S_1},\cdots,e^{S_I}\Sigma_Ie^{S_I}).$$ Then, for any convex sets $K_1,\cdots,K_I \subset \mathbb{R}^n$ such that $$\int_{K_i}x\, d\gamma_{\Sigma_i}(x)=0 \quad \forall 1\le i \le I,$$ it holds that $$\gamma_\Gamma \Big(\sum_{i=1}^I \lambda_i e^{S_i}K_i\Big) \ge \prod_{i=1}^I \gamma_{\Sigma_i}(K_i)^{\lambda_i}.$$
    Else, assume that \begin{equation} \label{Karcher_Condition_LpBM} \Gamma \le \prod_{i=1}^I \Big(\frac{\tilde \lambda_i}{\lambda_i}\Big)^{2(\frac1p-1)\tilde \lambda_i}\mathcal O_{\tilde \lambda}(e^{S_1}\Sigma_1e^{S_1},\cdots,e^{S_I}\Sigma_I e^{S_I}).\end{equation}
    Then, for any convex sets $K_1,\cdots,K_I \subset \mathbb{R}^n$  such that $$\int_{K_i}x\, d\gamma_{\Sigma_i}(x)=0 \quad \forall 1\le i \le I,$$ it holds that $$\gamma_\Gamma \Big(\sum_{i=1}^I \lambda_i e^{S_i}K_i\Big) \ge \prod_{i=1}^I \Big(\frac{\lambda_i}{\tilde \lambda_i}\Big)^{\frac{\tilde \lambda_i n}p} \prod_{i=1}^I \gamma_{\Sigma_i}(K_i)^{\tilde \lambda_i}.$$
Consequently, if \cref{Karcher_Condition_LpBM} holds for any $\tilde \lambda \in \Delta_I^\circ$, \begin{equation} \label{Geometric_BBL} \gamma_\Gamma\Big(\sum_{i=1}^I \lambda_i e^{S_i}K_i\Big)^{1/n} \ge \mathcal M_p^\lambda \big(\gamma_{\Sigma_1}(K_1)^{1/n},\cdots,\gamma_{\Sigma_I}(K_I)^{1/n}\big).\end{equation}
\end{theorem}

\hfill \\
\begin{remark} \label{Relaxed_Conditions_LpBM}
First, let us note that it would suffice to prove the statement for $S_i=0$ for any $1 \le i \le I.$ Indeed, it would be enough to apply such a statement to the covariances $\Gamma, e^{S_1}\Sigma_1e^{S_1},\cdots,e^{S_I}\Sigma_Ie^{S_I}$ and the convex sets $e^{S_1}K_1,\cdots,e^{S_I}K_I$ which satisfy $$\int_{e^{S_i}K_i}x\, d\gamma_{e^{S_i}\Sigma_ie^{S_i}}(x)=e^{S_i}\int_{K_i}y\, d\gamma_{\Sigma_i}(y)=0 \quad \forall 1 \le i \le I.$$ However, for convenience, we prefer to formulate this theorem in the most transparent way to better identify its strength. 

\hfill \\
Furthermore, a sufficient way of checking that \cref{Karcher_Condition_LpBM} holds for any $\tilde \lambda \in \Delta_I^\circ$ is to require $$\sum_{i=1}^I \lambda_i \Big(\Gamma^{-1/2}e^{S_i}\Sigma_i e^{S_i}\Gamma^{-1/2}\Big)^{q/2} \le I_n, \quad \text{ where } q=\frac p{p-1} \in (-\infty;0),$$ for $p \in (0;1).$ At $p=1$, a sufficient condition is $\Gamma^{-1/2}e^{S_i}\Sigma_i e^{S_i}\Gamma^{-1/2} \ge I_n$ for any $1 \le i \le I$, while at $p=0$, $$\sum_{i=1}^I \lambda_i \log\Big(\Gamma^{-1/2}e^{S_i}\Sigma_i e^{S_i}\Gamma^{-1/2} \Big) \ge 0.$$ 

\hfill \\
In fact, by \cref{Weighted_Positive_Sum}, $$\sum_{i=1}^I \lambda_i \log\Big(\Gamma^{-1/2}e^{S_i}\Sigma_ie^{S_i}\Gamma^{-1/2}\Big) \ge 0$$ implies that
$$\begin{array}{lll}
I_n &\le& \ds \mathcal O_\lambda(\Gamma^{-1/2}e^{S_1}\Sigma_1e^{S_1}\Gamma^{-1/2},\cdots,\Gamma^{-1/2}e^{S_I}\Sigma_Ie^{S_I}\Gamma^{-1/2}) \\ \\
\nm 
&\le& \ds \Gamma^{-1/2}\mathcal O_\lambda(e^{S_1}\Sigma_1e^{S_1},\cdots,e^{S_I}\Sigma_Ie^{S_I})\Gamma^{-1/2},
\end{array}$$ as required. At $p=1$, \cref{Karcher_Condition_LpBM} translates to $$\Gamma \le \mathcal O_{\tilde \lambda}(e^{S_1}\Sigma_1e^{S_1},\cdots,e^{S_I}\Sigma_Ie^{S_I}) \quad \forall \tilde \lambda \in \Delta_I^\circ,$$ which is equivalent to $\log(\Gamma^{-1/2}e^{S_i}\Sigma_ie^{S_i}\Gamma^{-1/2}) \ge 0$ for any $1 \le i \le I.$
Finally, for $p \in (0;1)$, if $$\sum_{i=1}^I \lambda_i (\Gamma^{-1/2}e^{S_i}\Sigma_ie^{S_i}\Gamma^{-1/2})^{q/2} \le I_n,$$ then $$\sum_{i=1}^I \tilde \lambda_i \Big( \Big(\frac{\lambda_i}{\tilde \lambda_i}\Big)^{\frac2q}\Gamma^{-1/2}e^{S_i}\Sigma_i e^{S_i}\Gamma^{-1/2} \Big)^{q/2} \le I_n \Longrightarrow \sum_{i=1}^I \tilde \lambda_i \log\Big[ \Big(\frac{\lambda_i}{\tilde \lambda_i}\Big)^{\frac2q}\Gamma^{-1/2}e^{S_i}\Sigma_i e^{S_i}\Gamma^{-1/2} \Big] \ge 0,$$ by concavity of the logarithm and as $q<0.$ Hence, by \cref{Weighted_Positive_Sum}, $$I_n \le \prod_{i=1}^I \Big(\frac{\tilde \lambda_i}{\lambda_i}\Big)^{2(\frac1p-1)\tilde \lambda_i}\mathcal O_{\tilde \lambda}(\Gamma^{-1/2}e^{S_1}\Sigma_1e^{S_1}\Gamma^{-1/2},\cdots,\Gamma^{-1/2}e^{S_I}\Sigma_Ie^{S_I}\Gamma^{-1/2}),$$ and this gives the desired condition. 
\end{remark}

\hfill \\
To consider inequalities with barycentric corrections, we introduce a new slicing approach in Section 4.1, which allows us to apply centered Gaussian inequalities in spaces of increasing dimensions and then pass to the limit. In the context of the Gaussian Brunn-Minkowski inequality, this gives the following extension of \cref{Han_Zhu_BM} to arbitrary finite families of sets with convexity only required for the weighted sum. 

\begin{theorem} \label{Generalized_Gaussian_BM}
Let $\lambda \in \Delta_I^\circ$ and $K_i \subset \R^n$ measurable sets of positive Gaussian measure such that $\sum_{i=1}^I \lambda_i K_i$ is convex. Then, $$e^{|\sum_{i=1}^I \lambda_i b_{K_i}|^2/2n} \gamma_n \Big(\sum_{i=1}^I \lambda_i K_i \Big)^{1/n} \ge \sum_{i=1}^I \lambda_i e^{|b_{K_i}|^2/2n} \gamma_n(K_i)^{1/n},$$ denoting by $$b_{K_i}:=\mathrm{bar}_\gamma(K_i)=\frac1{\gamma_n(K_i)}\int_{K_i}x\, d\gamma_n(x) \quad \forall 1 \le i \le I.$$ More generally, for $\Sigma \in \mathrm{Sym}_{>0}(n)$, $$e^{\langle \Sigma^{-1}\sum_{i=1}^I \lambda_i b_{K_i,\Sigma},\sum_{i=1}^I \lambda_i b_{K_i,\Sigma}\rangle/2n} \gamma_\Sigma \Big(\sum_{i=1}^I \lambda_i K_i\Big)^{1/n} \ge \sum_{i=1}^I \lambda_i e^{\langle \Sigma^{-1}b_{K_i,\Sigma},b_{K_i,\Sigma}\rangle/2n}\gamma_\Sigma(K_i)^{1/n},$$ where now $$b_{K_i,\Sigma}:=\mathrm{bar}_{\gamma_\Sigma}(K_i)=\frac1{\gamma_\Sigma(K_i)}\int_{K_i}x\, d\gamma_{\Sigma}(x) \quad \forall 1 \le i \le I.$$ 
\end{theorem}
\hfill \\
The case $I=2$ is exactly \cite[Corollary 1.6]{Han_Zhu} when assuming $K_1$ and $K_2$ are convex (see \cref{Han_Zhu_BM}). For $I \ge 3$, the inequality does not follow by direct iteration of the statement for two sets as the Gaussian barycenter of a Minkowski sum is not determined by the barycenters of its summands. In particular, it is new to the author's knowledge.

\hfill \\
Note that Han-Zhu initially study Poincar\'e estimates of Wasserstein geodesics between probability measures satisfying strong log-concavity in \cite{Han_Zhu} (see \cite{Gozlan_Malamut_Waldspurger} for an extension of this result to the sharp Talagrand inequality). In particular, a barycentric extension of these results (see \cite[Question 1.7]{Han_Zhu}) would be of great interest and if true, would recover \cref{Generalized_Gaussian_BM} using the slicing method introduced here. 

\subsection{Structure of the paper}

The second section is devoted to the proof of \cref{Sharp_Constant_BM_Ineq}, which is based on the Gaussian saturation principle of the associated functional form, which is a Gaussian Forward-Reverse Brascamp-Lieb. 
In Section 3, we provide functional extensions and study equality cases first when $p=0$. Based on this case, we conclude the proof of \cref{Main_Thm_IneqPart} and settle the equality cases of the Borell-Brascamp-Lieb inequality corresponding to \cref{Geometric_BBL}. In addition, we answer a question formulated by Faifman-Polterovich (see \cite[Question 5.3]{Faifman_Polt}) in the negative. This question asks whether the scalar (B)-property, or more generally its symmetric-matrix version, holds for convex bodies with Lebesgue barycenters at the origin. In the last section, our new slicing method is introduced, which enables us to treat the Gaussian Brunn-Minkowski for not-necessarily centered sets and provide another application through an improved Gaussian correlation inequality. Thanks to further variational arguments, we give a complete rigidity statement for the two considered inequalities mentioned above.

\section{Proof of \cref{Sharp_Constant_BM_Ineq}}

The main functional tool used in this paper is the Gaussian Forward-Reverse Brascamp-Lieb of \cite{Strong_Gaussian_Correlation}. 

\begin{theorem}
    Let $L_j: \oplus_{i=1}^I \mathbb{R}^{n_i} \rightarrow \mathbb{R}^{m_j}$ be a surjective linear map for any $1 \le j \le J$, $c_i>0$ for any $1 \le i \le I$ and $d_j>0$ for any $1 \le j \le J.$ Further, consider $\Sigma_i \in \mathrm{Sym}_{>0}(n_i), \Gamma_j \in \mathrm{Sym}_{\ge 0}(m_j)$ for any $1\le i\le I, 1 \le j \le J$ and $\mathcal Q \in \mathrm{Sym}(\oplus_{i=1}^I \mathbb{R}^{n_i})$. Then, for all log-concave functions $f_i$ that are centered with respect to $\gamma_{\Sigma_i}$, \textit{i.e.} $$\int_{\mathbb{R}^{n_i}} xf_i(x)\, d\gamma_{\Sigma_i}(x)=0 \quad \forall 1 \le i \le I,$$ and log-concave functions $h_j$, with positive finite Gaussian integrals, $$\prod_{i=1}^I f_i(x_i)^{c_i} \le e^{-\frac12 \langle \mathcal Qx,x\rangle}\prod_{j=1}^J h_j(L_jx)^{d_j} \Longrightarrow \prod_{i=1}^I \Big(\int_{\mathbb{R}^{n_i}}f_i\, d\gamma_{\Sigma_i}\Big)^{c_i} \le C_{\mathcal G}\prod_{j=1}^J \Big(\int_{\mathbb{R}^{m_j}}h_j\, d\gamma_{\Gamma_j}\Big)^{d_j},$$ where $$\begin{array}{lll}
    C_{\mathcal G} &=& \ds \sup_{\substack{A_i,B_j \ge 0 \\ \mathrm{diag}(c_iA_i) \ge \mathcal Q+\sum_{j=1}^J d_jL_j^*B_jL_j}} \frac{\prod_{j=1}^J \det(I_{m_j}+\Gamma_j B_j)^{d_j/2}}{\prod_{i=1}^I \det(I_{n_i}+\Sigma_i A_i)^{c_i/2}}
    \\ \\
    \nm 
    &=& \ds \sup_{\substack{\prod_{i=1}^I g_{A_i}^{c_i}(x_i) \le e^{-\frac12 \langle \mathcal Qx,x\rangle} \prod_{j=1}^J g_{B_j}(L_jx)^{d_j} \\ A_i,B_j \ge 0}} \frac{\prod_{i=1}^I (\int_{\mathbb{R}^{n_i}}g_{A_i}\,d\gamma_{\Sigma_i})^{c_i}}{\prod_{j=1}^J (\int_{\mathbb{R}^{m_j}}g_{B_j}\, d\gamma_{\Gamma_j})^{d_j}}
    \end{array}$$ is the best constant, when $f_i$ and $h_j$ are restricted to centered Gaussians in the implication above, defining $$g_A(x)=e^{-\frac12 \langle Ax,x\rangle}.$$  
\end{theorem}

\begin{proof}[Proof of \cref{Sharp_Constant_BM_Ineq}]
    The goal is to compute the Gaussian constant of the following Forward-Reverse Brascamp-Lieb: for $f_i$ log-concave and centered with respect to $\gamma_{\Sigma_i}$ for $ 1 \le i \le I$ and $h$ log-concave, $$\prod_{i=1}^I f_i(x_i)^{\tilde \lambda_i} \le h\Big(\sum_{i=1}^I \lambda_i x_i\Big) \Longrightarrow \prod_{i=1}^I \Big(\int_{\mathbb{R}^n}f_i\, d\gamma_{\Sigma_i}\Big)^{\tilde \lambda_i} \le C_{\mathrm{FRBL}}\int_{\mathbb{R}^n}h \, d\gamma_\Gamma.$$ Then, $$C_{\mathrm{FRBL}}=\sup \Bigg\{\frac{\det(I_n+\Gamma B)}{\prod_{i=1}^I \det(I_n+\Sigma_i A_i)^{\tilde \lambda_i}} : \mathrm{diag}(\tilde \lambda_1 A_1,\cdots,\tilde \lambda_I A_I) \ge (\lambda_i \lambda_j B)_{1 \le i,j \le I} \Bigg\}^{1/2}.$$ For fixed $A_i>0$, the largest possible $B$ is $$\Big(\sum_{i=1}^I \frac{\lambda_i^2}{\tilde \lambda_i}A_i^{-1}\Big)^{-1},$$ as it is easy to show that $$\min_{\sum_{i=1}^I \lambda_ix_i=y}\sum_{i=1}^I \tilde \lambda_i \langle A_ix_i,x_i\rangle=\Big\langle\Big(\sum_{i=1}^I \frac{\lambda_i^2}{\tilde \lambda_i}A_i^{-1}\Big)^{-1}y,y\Big\rangle.$$ Consequently, setting $$A_i=\Big(\frac{\lambda_i}{\tilde \lambda_i}\Big)^2X_i, \quad Y_i=\Big(\frac{\lambda_i}{\tilde \lambda_i}\Big)^2 \Sigma_i \text{ and } H_{\tilde \lambda}(X_1,\cdots,X_I)=\Big(\sum_{i=1}^I \tilde \lambda_i X_i^{-1}\Big)^{-1}$$ yields that $$C_{\mathrm{FRBL}}^2=\sup_{X_1,\cdots,X_I>0} \frac{\det(I_n+\Gamma H_{\tilde \lambda}(X_1,\cdots,X_I))}{\prod_{i=1}^I \det(I_n+Y_iX_i)^{\tilde \lambda_i}}.$$ 
    We claim that this is equal to $$\sup_{X>0}\frac{\det(I_n+\Gamma X)}{\prod_{i=1}^I \det(I_n+Y_iX)^{\tilde \lambda_i}}=\max_{E \text{ linear subspace of } \mathbb{R}^n}\frac{\det(\Gamma_E)}{\prod_{i=1}^I (\det(Y_i)_E)^{\tilde \lambda_i}}.$$ 
    For a decomposition $\mathbb{R}^n=E \oplus F$, we define the Schur complement $$S_E(A):=A_F-A_{FE}A_E^{-1}A_{EF},$$ where we write $$A=\begin{pmatrix}
        A_E & A_{EF} \\
        A_{FE} & A_F
    \end{pmatrix}.$$ Its well-known variational characterization is given by $$\langle S_E(A)v,v\rangle=\inf_{u \in E}\langle A(u+v),u+v\rangle, \quad v \in F.$$ It follows that $S_E$ is monotone and concave, \textit{i.e.} $$A \le B \Longrightarrow S_E(A) \le S_E(B), \quad S_E\Big(\sum_{i=1}^I \lambda_i A_i\Big) \ge \sum_{i=1}^I \lambda_i S_E(A_i).$$ We will use the following two inequalities: first, \begin{equation} \label{Proj_Quotient_Det_Decreases} \frac{\det(A+Y)}{\det A} \ge \frac{\det(A_E+Y_E)}{\det A_E} \quad \forall A,Y>0,\end{equation} since $$\frac{\det(A+Y)}{\det A}=\frac{\det(A_E+Y_E)}{\det A_E} \cdot \frac{\det S_E(A+Y)}{\det S_E(A)},$$ and one may use monotonicity of $S_E.$ Second, \begin{equation} \label{Proj_Convex_Sum_Det_Increases} \frac{\prod_{i=1}^I (\det A_i)^{\tilde \lambda_i}}{\det\Big(\sum_{i=1}^I \tilde \lambda_i A_i\Big)} \le \frac{\prod_{i=1}^I (\det (A_i)_F)^{\tilde \lambda_i}}{\det\Big(\sum_{i=1}^I \tilde \lambda_i (A_i)_F\Big)}.\end{equation} To see this, we may use $$\frac{\prod_{i=1}^I (\det A_i)^{\tilde \lambda_i}}{\det\Big(\sum_{i=1}^I \tilde \lambda_i A_i\Big)}=\frac{\prod_{i=1}^I (\det (A_i)_F)^{\tilde \lambda_i}}{\det\Big(\sum_{i=1}^I \tilde \lambda_i (A_i)_F\Big)} \cdot \frac{\prod_{i=1}^I \det(S_F(A_i))^{\tilde \lambda_i}}{\det S_F\Big(\sum_{i=1}^I \tilde \lambda_i A_i\Big)},$$ and by concavity of $\log \det$, $$\det S_F\Big(\sum_{i=1}^I \tilde \lambda_i A_i\Big) \ge \det \Big( \sum_{i=1}^I \tilde \lambda_i S_F(A_i) \Big) \ge \prod_{i=1}^I \det(S_F(A_i))^{\tilde \lambda_i}.$$ 
    \hfill \\
    Setting $$Q_i=X_i^{-1}, \quad \overline{Q}=\sum_{i=1}^I \tilde \lambda_i Q_i,$$ we have that $$\frac{\det(I_n+\Gamma H_{\tilde \lambda}(X_1,\cdots,X_I))}{\prod_{i=1}^I \det(I_n+Y_iX_i)^{\tilde \lambda_i}}=\frac{\det(\overline Q+\Gamma)}{\det(\overline Q)} \prod_{i=1}^I \Big(\frac{\det Q_i}{\det(Q_i+Y_i)}\Big)^{\tilde \lambda_i}:=R(Q_1,\cdots,Q_I).$$  
    \hfill \\
    Observe that every interior critical point of $R$ has value at most $1.$ Indeed, differentiating $\log R$ with respect to each $Q_i$ gives $$D:=\overline{Q}^{-1}-(\overline Q+\Gamma)^{-1}=Q_i^{-1}-(Q_i+Y_i)^{-1} \quad \forall 1 \le i \le I.$$ Then, $\Gamma>0$ implies that $D>0$ and $$0<D^{1/2}Q_iD^{1/2}<I_n, \quad \frac{\det(Q_i+Y_i)}{\det(Q_i)}=\det(I_n-D^{1/2}Q_iD^{1/2})^{-1},$$ while $$0<D^{1/2}\overline Q D^{1/2}<I_n, \quad \frac{\det(\overline Q+\Gamma)}{\det(\overline Q)}=\det(I_n-D^{1/2}\overline Q D^{1/2})^{-1}.$$ Therefore, $$R(Q_1,\cdots,Q_I)=\frac{\prod_{i=1}^I \det(I_n-D^{1/2}Q_iD^{1/2})^{\tilde \lambda_i}}{\det(I_n-D^{1/2}\overline QD^{1/2})} \le 1,$$ by concavity of $\log \det$ and $$I_n-D^{1/2}\overline QD^{1/2}=\sum_{i=1}^I \tilde \lambda_i(I_n-D^{1/2}Q_iD^{1/2}).$$ 
    In particular, it suffices to control the values of $R$ over the boundary of the positive definite cone. 
    \hfill \\
    Define $$M:=\max_{E \text{ linear subspace of } \mathbb{R}^n} \frac{\det(\Gamma_E)}{\prod_{i=1}^I (\det(Y_i)_E)^{\tilde \lambda_i}} \ge 1,$$ by considering the zero subspace. We prove that $$S=\sup_{Q_i>0} R(Q_1,\cdots,Q_I) \le M$$ by induction on $n.$ Dimension $0$ is trivial, so let us assume the assertion in smaller dimension.
    Choose a maximizing sequence $Q_i^{(k)}$ and note that $S \ge 1$ by considering $(tI_n)_{1 \le i \le I}$ and letting $t \rightarrow +\infty.$ 
    After passing to a subsequence, every ordered eigenvalue of $$\overline{Q}^{(k)}=\sum_{i=1}^I \tilde \lambda_i Q_i^{(k)}$$ has a limit in $[0;+\infty].$ 
    \hfill \\
    We split the analysis in three cases. 
    First, assume that some eigenvalues of $\overline{Q}^{(k)}$ diverge. Let $E_k$ be the spectral subspace corresponding to the eigenvalues that remain bounded and $F_k=E_k^\perp.$ Then, $\dim(E_k)=r<n$, and in this decomposition, $$\overline{Q}^{(k)}=U_k \oplus V_k, \quad \lambda_{\min}(V_k)=\min_{\|x\|=1, x \in F_k}\langle V_kx,x\rangle \rightarrow_{k \rightarrow +\infty}+\infty.$$ 
    \hfill \\
    Then, $$\frac{\det(\overline Q^{(k)}+\Gamma)}{\det \overline Q^{(k)}}=\frac{\det(U_k+\Gamma_{E_k})}{\det U_k}c_k,$$ where $$c_k=\frac{\det(V_k+\Gamma_{F_k}-\Gamma_{F_kE_k}(U_k+\Gamma_{E_k})^{-1}\Gamma_{E_kF_k})}{\det V_k}.$$ As the matrix added to $V_k$ is positive semi-definite and less than $\Gamma_{F_k}$, $$1 \le c_k \le \Big(1+\frac{\|\Gamma\|_{\mathrm{op}}}{\lambda_{\min}(V_k)}\Big)^{n-r} \rightarrow_{k \rightarrow +\infty}1.$$ By \cref{Proj_Quotient_Det_Decreases}, $$\frac{\det(Q_i^{(k)}+Y_i)}{\det Q_i^{(k)}} \ge \frac{\det\big((Q_i^{(k)})_{E_k}+(Y_i)_{E_k} \big)}{\det \big((Q_i^{(k)})_{E_k}\big)},$$ which implies that $$R(Q_1^{(k)},\cdots,Q_I^{(k)}) \le c_kR_{E_k}((Q_i^{(k)})_{E_k}).$$ We may apply the induction hypothesis to $E_k$, so that  $R_{E_k}((Q_i^{(k)})_{E_k})$ is bounded above by $$\max_{V \text{ linear subspace of } E_k}\frac{\det \Gamma_V}{\prod_{i=1}^I (\det(Y_i)_V)^{\tilde \lambda_i}} \le M.$$ Thus, $$\limsup_{k \rightarrow +\infty}R(Q_1^{(k)},\cdots,Q_I^{(k)}) \le M.$$ 
    \hfill \\
    Now, assume that $\overline{Q}^{(k)}$ is bounded and has a positive definite limit. As $$0<Q_i^{(k)} \le \tilde \lambda_i^{-1}\overline Q^{(k)},$$ all the matrices are bounded and up to a further subsequence, $$Q_i^{(k)} \rightarrow Q_i^{(0)} \ge 0, \quad \overline{Q}^{(k)} \rightarrow \overline Q^{(0)}>0 \text{ as } k \rightarrow +\infty.$$ If any $Q_i^{(0)}$ is singular, $R(Q_1^{(k)},\cdots,Q_I^{(k)}) \rightarrow_{k \rightarrow +\infty}0.$ Otherwise, all $Q_i^{(0)}>0$, so that continuity shows that the supremum is attained at an interior tuple $(Q_i)_{1 \le i \le I}.$ By what we previously showed, $$S \le 1 \le M.$$ 
    \hfill \\
    Finally, assume that $\overline{Q}^{(k)}$ is bounded and has a singular limit. Again, up to a subsequence, $Q_i^{(k)} \rightarrow Q_i^{(0)} \ge 0$, and set $$E=\ker \overline Q^{(0)}, \quad F=E^\perp.$$ Then, every $Q_i^{(0)}$ vanishes on $E$, since $\overline Q^{(0)}=\sum_{i=1}^I \tilde \lambda_i Q_i^{(0)}.$ Hence, $$Q_i^{(0)}=0_E \oplus Q_i', \quad \overline Q^{(0)}=0_E \oplus \overline Q',$$ with $Q_i' \ge 0$ and  $\overline Q'=\sum_{i=1}^I \tilde \lambda_i Q_i'>0$ on $F.$ Hence, by \cref{Proj_Convex_Sum_Det_Increases}, $$\begin{array}{lll}
    R(Q_1^{(k)},\cdots,Q_I^{(k)}) &=& \ds \frac{\det(\overline Q^{(k)}+\Gamma)}{\prod_{i=1}^I \det(Q_i^{(k)}+Y_i)^{\tilde \lambda_i}} \cdot \frac{\prod_{i=1}^I \det(Q_i^{(k)})^{\tilde \lambda_i}}{\det(\overline Q^{(k)})}
    \\ \\
    \nm 
    &\le& \ds  \frac{\det(\overline Q^{(k)}+\Gamma)}{\prod_{i=1}^I \det(Q_i^{(k)}+Y_i)^{\tilde \lambda_i}} \cdot  \frac{\prod_{i=1}^I \det((Q_i^{(k)})_F)^{\tilde \lambda_i}}{\det((\overline Q^{(k)})_F)}.
    \end{array}$$ If some $Q_i'$ is singular, the right-hand side tends to zero and otherwise all $Q_i'>0$ so that we may pass to the limit. 
    Hence, $$\limsup_{k \rightarrow +\infty} R(Q_i^{(k)}) \le \frac{\det \Gamma_E}{\prod_{i=1}^I (\det(Y_i)_E)^{\tilde \lambda_i}} R_F(Q_i';S_E(\Gamma),S_E(Y_i)).$$ The induction hypothesis applies on $F$, whose dimension is smaller than $n.$ Consequently, the right-hand side is at most $$\max_{V \text{ linear subspace of } F} \frac{\det \Gamma_{E \oplus V}}{\prod_{i=1}^I (\det(Y_i)_{E \oplus V})^{\tilde \lambda_i}} \le M,$$ using that $$\det \Gamma_E \det \Big((S_E(\Gamma))_V\Big)=\det \Gamma_{E \oplus V},$$ and similarly for $Y_i.$  
    \hfill \\
    We conclude that $S \le M.$ 
    If we set $$X_1=\cdots=X_I=X_t:=t P_E+t^{-1}P_{E^\perp}, \quad t \rightarrow +\infty,$$ then $$H_{\tilde \lambda}(X_t,\cdots,X_t)=X_t,$$ and $$\det(I_n+AX_t) \sim t^{\dim(E)}\det(A_E) \quad \forall A \in \mathrm{Sym}_{>0}(n).$$ It follows that $$\frac{\det(I_n+\Gamma X_t)}{\prod_{i=1}^I \det(I_n+Y_iX_t)^{\tilde \lambda_i}} \rightarrow_{t \rightarrow +\infty} \frac{\det \Gamma_E}{\prod_{i=1}^I (\det (Y_i)_E)^{\tilde \lambda_i}}.$$  
    \hfill \\
    This proves the desired claim and recalling that $Y_i=\Big(\frac{\lambda_i}{\tilde \lambda_i}\Big)^2 \Sigma_i$, concludes the proof of the theorem. 
\end{proof}

\section{Interpolation between (B)-theorem and Gaussian Brunn-Minkowski}

In this section, we study in detail the Gaussian inequalities obtained in the previous section, that have sharp constant exactly equal to $1$ under a natural matrix condition satisfied by the covariances. 

\subsection{The bilinear (B)-theorem and study of equality cases}

\hfill \\
Let us first consider a version of the bilinear (B)-theorem with multiple sets, \textit{i.e.} the part $p=0$ of \cref{Main_Thm_IneqPart}. Although the next theorem is in fact contained in \cref{Sharp_Constant_BM_Ineq}, we will provide a self-contained computation of the Gaussian constant to be able to treat the equality cases. 

\begin{theorem} \label{Multiple_BThm}
    If $I \ge 2, \lambda \in \Delta_I^\circ$ and $S_1,\cdots,S_I \in \mathrm{Sym}(n)$ such that $\sum_{i=1}^I \lambda_i S_i \ge 0.$ Then, for $K_1,\cdots,K_I$ convex and centered with respect to the standard Gaussian measure, $$\gamma_n \Big(\sum_{i=1}^I \lambda_i e^{S_i}K_i) \ge \prod_{i=1}^I \gamma_n(K_i)^{\lambda_i}.$$ More generally, if $\Gamma, \Sigma_i \in \mathrm{Sym}_{>0}(n)$ such that $$\Gamma \le \mathcal O_\lambda(e^{S_1}\Sigma_1e^{S_1},\cdots,e^{S_I}\Sigma_Ie^{S_I}),$$ and $K_1,\cdots,K_I$ convex and with $K_i$ being centered with respect to $\gamma_{\Sigma_i}$ for any $1 \le i \le I$, then $$\gamma_\Gamma\Big(\sum_{i=1}^I \lambda_i e^{S_i}K_i\Big) \ge \prod_{i=1}^I \gamma_{\Sigma_i}(K_i)^{\lambda_i}.$$
\end{theorem}

\begin{remark}
    The ``More generally'' is because if we consider $\Gamma=\Sigma_i=I_n$, \begin{equation} \label{Weighted_Positive_Sum} \sum_{i=1}^I \lambda_i S_i \ge 0 \Longrightarrow  I_n \le \mathcal O_\lambda(e^{2S_1},\cdots,e^{2S_I}).\end{equation}
    Indeed, assume that $\sum_{i=1}^I \lambda_i \log(A_i) \ge 0$ and let $M=\mathcal O_\lambda(A_1,\cdots,A_I).$ If $\mu$ is the smallest eigenvalue of $M$ and $v \ne 0$ is a corresponding unit vector, $$\begin{array}{lll}
    \ds \big \langle v,\log(M^{-1/2}A_iM^{-1/2})v\big\rangle &=& \ds \int_0^{+\infty}\Big(\frac1{1+t}-\mu\langle (A_i+tM)^{-1}v,v\rangle \Big)\, dt \\ \\
    \nm 
    &\ge&\ds \langle \log(A_i)v,v\rangle-\log \mu,
    \end{array}$$ using here that $M \ge \mu I_n$ and for $x>0$, $$\log(x)=\int_0^{+\infty}\Big( \frac1{1+t}-\frac1{x+t}\Big)\, dt.$$
    Then, taking the weighted sum shows that $$0 \ge \Big\langle \sum_{i=1}^I \lambda_i \log(A_i)v,v\Big\rangle-\log(\mu) \ge -\log(\mu) \Longrightarrow \mu \ge 1.$$ 
\end{remark}

\begin{proof}
    This is a consequence of the following forward-reverse Brascamp-Lieb (see \cite[Theorem 1.10]{Strong_Gaussian_Correlation}), with the map $$L:(\mathbb{R}^n)^I \rightarrow \mathbb{R}^n, \quad (x_1,\cdots,x_I) \mapsto \sum_{i=1}^I \lambda_i e^{S_i}x_i \text{ being surjective.}$$ This implies that if $f_i$ is log-concave, centered with respect to $\gamma_{\Sigma_i}$ for any $1 \le i \le I$ and $h$ is log-concave, \begin{equation} \label{FRBL_MultipleBthm} \prod_{i=1}^I f_i^{\lambda_i}(x_i) \le h \Big(\sum_{i=1}^I \lambda_i e^{S_i}x_i\Big) \Longrightarrow \prod_{i=1}^I \Big(\int_{\R^n}f_i\, d\gamma_{\Sigma_i} \Big)^{\lambda_i} \le \int_{\R^n}h\, d\gamma_\Gamma.\end{equation} Indeed, the Gaussian constant is equal to $$\sup \left\{ \frac{\det(I_n+\Gamma B)}{\prod_{i=1}^I \det(I_n+\Sigma_i A_i)^{\lambda_i}} : \mathrm{diag}(\lambda_1 A_1,\cdots,\lambda_I A_I) \ge (\lambda_i \lambda_j e^{S_i}Be^{S_j})_{1 \le i,j \le I}, A_i>0,B>0 \right\}^{1/2}=1,$$ by Proposition \ref{Gaussian_Case_MultipleBThm}. 
\end{proof}

\begin{remark}
    As already explained previously, we may deduce \cref{Multiple_BThm} from \cref{Sharp_Constant_BM_Ineq}. Indeed, for every subspace $E \subset \mathbb{R}^n$ and $\tilde \lambda \in \Delta_I^\circ$, $$\mathcal O_{\tilde \lambda}(\Sigma_1,\cdots,\Sigma_I)_E \le \mathcal O_{\tilde \lambda}((\Sigma_1)_E,\cdots,(\Sigma_I)_E).$$
    This holds by considering $U=I_E \oplus (-I_F)$ and noticing that $$\frac{\Sigma_i+U\Sigma_iU}{2}=(\Sigma_i)_E \oplus (\Sigma_i)_F,$$ together with the joint concavity of the Karcher mean $$\mathcal O_{\tilde \lambda}\Big(\Big(\frac{\Sigma_i+U\Sigma_iU}{2}\Big)_{1 \le i \le I} \Big) \ge \frac{\mathcal O_{\tilde \lambda}((\Sigma_i)_{1 \le i \le I})+U\mathcal O_{\tilde \lambda}((\Sigma_i)_{1 \le i \le I})U}{2}.$$ Combining these observations and using that the Karcher mean preserves the block structure, $$\mathcal O_{\tilde \lambda}\big[ ((\Sigma_i)_E)_{1 \le i \le I} \big] \oplus \mathcal O_{\tilde \lambda}\big[ ((\Sigma_i)_F)_{1 \le i \le I} \big] \ge \mathcal O_{\tilde \lambda}\big[(\Sigma_i)_{1 \le i \le I} \big]_E \oplus \mathcal O_{\tilde \lambda}\big[(\Sigma_i)_{1 \le i \le I} \big]_F.$$ 
    \hfill \\
    Consequently, if $\Gamma \le \mathcal O_{\tilde \lambda}(\Sigma_1,\cdots,\Sigma_I)$, $$\Gamma_E \le \mathcal O_{\tilde \lambda}(\Sigma_1,\cdots,\Sigma_I)_E \le \mathcal O_{\tilde \lambda}((\Sigma_1)_E,\cdots,(\Sigma_I)_E),$$ so that taking determinants, we have $$\det(\Gamma_E) \le \prod_{i=1}^I (\det((\Sigma_i)_E))^{\tilde \lambda_i},$$ by the determinant identity for the Karcher mean. Thus, the constant of \cref{Sharp_Constant_BM_Ineq} is less or equal in this case to $$\max_{E \text{ linear subspace of } \mathbb{R}^n} \prod_{i=1}^I \Big(\frac{\tilde \lambda_i}{\lambda_i}\Big)^{\dim(E)\tilde \lambda_i}=\prod_{i=1}^I\Big(\frac{\tilde \lambda_i}{\lambda_i}\Big)^{n\tilde \lambda_i}.$$ In particular, this is equal to $1$ if $\tilde \lambda=\lambda$ and implies that if $\Gamma \le \mathcal O_\lambda(\Sigma_1,\cdots,\Sigma_I)$, $$\gamma_\Gamma\Big(\sum_{i=1}^I\lambda_i K_i\Big) \ge \prod_{i=1}^I \gamma_{\Sigma_i}(K_i)^{\lambda_i}.$$ It then suffices to apply this inequality to the convex set $e^{S_i}K_i$ that is centered with respect to $\gamma_{e^{S_i}\Sigma_ie^{S_i}}$ for any $1 \le i \le I$, as already explained in the introduction. 
\end{remark}

\hfill \\
To compute the Gaussian constant associated to the previous Forward-Reverse Brascamp-Lieb formulation, we first establish the following auxiliary lemma. 

\begin{lemma} \label{Aux_Lemma}
    If $C \in \R^{m \times n}, P_i\in \mathrm{Sym}(m)$ and $Q_i \in \mathrm{Sym}(n)$ for any $1 \le i \le I$, such that $$\sum_{i=1}^I \lambda_i P_i \ge 0, \sum_{i=1}^I \lambda_i Q_i \ge 0,$$ then $$\det(I_n+C^*C) \le \prod_{i=1}^I \det\Big( I_n+(e^{P_i}Ce^{Q_i})^* (e^{P_i}Ce^{Q_i}) \Big)^{\lambda_i}.$$
\end{lemma}

\begin{proof}
    For $P \in \mathrm{Sym}(m), Q \in \mathrm{Sym}(n)$, define $$C_t=e^{tP}Ce^{tQ}, \; F(t)=\log \det(I_n+C_t^*C_t) \quad \forall t \in \R.$$ Define $$Z=\begin{pmatrix}
        C \\
        I_n 
    \end{pmatrix}, R=\begin{pmatrix}
        P & 0_{m \times n} \\
        0_{n \times m} & -Q
    \end{pmatrix} \text{ and } Z_t=e^{tR}Z=\begin{pmatrix}
        e^{tP}C \\
        e^{-tQ}
    \end{pmatrix} \quad \forall t \in \R.$$ Then, $Z_t^*Z_t=C^*e^{2tP}C+e^{-2tQ}$, so that $$e^{tQ}Z_t^*Z_te^{tQ}=e^{tQ}C^*e^{tP}(e^{tP}Ce^{tQ})+I_n=C_t^*C_t+I_n.$$ Hence, $$F(t)=2\mathrm{Tr}(Q)t+\log \det(Z_t^*Z_t) \quad \forall t \in \R.$$
    We claim that $F$ is convex. Recall that $$\begin{array}{lll}
    \displaystyle \frac d{dt}\log \det(Z_t^*Z_t)=\mathrm{Tr} \Big( (Z_t^*Z_t)^{-1} \frac d{dt}(Z_t^*Z_t) \Big),\; \frac{d^2}{dt^2} \log \det(Z_t^*Z_t)=\mathrm{Tr}\Big( (Z_t^*Z_t)^{-1} \\ \\
    \displaystyle \frac{d^2}{dt^2}(Z_t^*Z_t) -\big [(Z_t^*Z_t)^{-1}\frac d{dt}(Z_t^*Z_t) \big ]^2 \Big).
    \end{array}$$ Now, $$\begin{array}{lll}
    \displaystyle \dot Z_t=RZ_t \Longrightarrow \frac d{dt}(Z_t^*Z_t)=(\dot Z_t)^* Z_t+Z_t^*\dot Z_t=2Z_t^*RZ_t, \\ \\
    \displaystyle \frac{d^2}{dt^2}(Z_t^*Z_t)=4Z_t^*R^2Z_t.
    \end{array}$$ Hence, defining $\Pi_t=Z_t(Z_t^*Z_t)^{-1}Z_t^*$, we notice that $$\Pi_t^2=Z_t(Z_t^*Z_t)^{-1}Z_t^*Z_t(Z_t^*Z_t)^{-1}Z_t^*=\Pi_t,$$ and therefore, $$\begin{array}{lll}
    \displaystyle \frac{d^2}{dt^2} \log \det(Z_t^*Z_t) &=& \displaystyle \mathrm{Tr}(4 \Pi_t R^2-4\Pi_t R \Pi_t R)=4 \mathrm{Tr}(\Pi_t R(I_{m+n}-\Pi_t)R) \\ \\ 
    &=& \displaystyle 4\|\Pi_t R(I_{m+n}-\Pi_t)\|_{\mathrm{HS}}^2 \ge 0,
    \end{array}$$ proving convexity. Moreover, $$\begin{array}{lll}
    \displaystyle F'(0)=2\mathrm{Tr}(Q)+2\mathrm{Tr}\Big( (I+C^*C)^{-1}(C^*PC-Q) \Big)=2\mathrm{Tr}(PU)+2\mathrm{Tr}(QV),
    \end{array}$$ where $$U=C(I_n+C^*C)^{-1}C^* \ge 0, V=I_n-(I_n+C^*C)^{-1} \ge 0.$$
    \hfill \\
    Now, define $$F_i(t)=\log \det \Big(I_n+(e^{tP_i}Ce^{tQ_i})^* (e^{tP_i}Ce^{tQ_i}) \Big) \quad \forall t \in \R, 1 \le i \le I.$$ 
    Then, by convexity of $F_i$, $$F_i(1) \ge F_i(0)+F_i'(0) \quad \forall 1 \le i \le I \;  \Longrightarrow \sum_{i=1}^I \lambda_i F_i(1) \ge \sum_{i=1}^I \lambda_i F_i(0)+\sum_{i=1}^I \lambda_i F_i'(0).$$ But, $F_i(0)=\log \det(I_n+C^*C)$ and $$\sum_{i=1}^I \lambda_i F_i'(0)=2\mathrm{Tr}\Big(\sum_{i=1}^I \lambda_i P_i U \Big)+2\mathrm{Tr}\Big(\sum_{i=1}^I \lambda_i Q_i V \Big) \ge 0,$$ as the trace of the product of two positive semi-definite matrices is always non-negative. Hence, $$\sum_{i=1}^I \lambda_i \log \det \Big(I_n+(e^{P_i}Ce^{Q_i})^* (e^{P_i}Ce^{Q_i})\Big) \ge \log \det(I_n+C^*C).$$ 
\end{proof}

\begin{remark}
    There is equality in Lemma \ref{Aux_Lemma} for $C$ with full column-rank if and only if $$\Big(\sum_{i=1}^I \lambda_i P_i\Big)C=0, \sum_{i=1}^I \lambda_i Q_i=0 \text{ and } P_iC+CQ_i=0 \quad \forall 1 \le i \le I.$$ Indeed, if there is equality, as $$F_i(1)-F_i(0)-F_i'(0)=\int_0^1 (1-t) F_i''(t)\, dt \ge 0,$$ by convexity, we get that $$\sum_{i=1}^I \lambda_i (F_i(1)-F_i(0))=2\mathrm{Tr}\Big(\sum_{i=1}^I \lambda_i P_iU \Big)+2\mathrm{Tr}\Big(\sum_{i=1}^I \lambda_i Q_iV\Big)+\sum_{i=1}^I \lambda_i \int_0^1 (1-t)F_i''(t)\, dt,$$ where $U=C(I_n+C^*C)^{-1}C^* \ge 0, V=I_n-(I_n+C^*C)^{-1}>0$, as $C$ is full-rank. So each term on the right-hand side must vanish, which implies that $$\sum_{i=1}^I \lambda_i Q_i=0, \text{ as } V>0,$$ and $$\sum_{i=1}^I \lambda_i P_i C=0,$$ as $U$ is positive definite on $\mathrm{Ran}(C).$ For the last term to be zero, we must in particular have that $$4\|\Pi_{0} R_i(I_{m+n}-\Pi_{0})\|_{\mathrm{HS}}^2=F_i''(0)=0,$$ where, for any $1 \le i \le I$, $$Z=\begin{pmatrix}
        C \\
        I_n
    \end{pmatrix},R_i=\begin{pmatrix}
        P_i & 0_{m \times n} \\
        0_{n \times m} & -Q_i
    \end{pmatrix} \text{ and } \Pi_0=Z(Z^*Z)^{-1}Z^*,$$ which is the orthogonal projection onto the range of $Z.$ This means that $\mathrm{Ran}(Z)$ stays invariant under $R_i$, so that for any $x \in \R^n$, there exists $y \in \R^n$ so that $$\begin{pmatrix}
        P_iCx \\
        -Q_ix
    \end{pmatrix}=R_i\begin{pmatrix}
        Cx \\
        x
    \end{pmatrix}=\begin{pmatrix}
        Cy \\
        y
    \end{pmatrix},$$ so that $P_iCx=-CQ_ix$ for any $x \in \R^n$, as desired.
\end{remark}
\hfill \\
We now compute the Gaussian constant and record the Gaussian equality cases, to later be able to study the equality cases of Theorem \ref{Multiple_BThm}. 

\begin{proposition} \label{Gaussian_Case_MultipleBThm}
    Fix $A_1,\cdots,A_I, B \in \mathrm{Sym}_{\ge 0}(n)$ such that $$\mathrm{diag}(\lambda_1A_1,\cdots,\lambda_I A_I) \ge (\lambda_i\lambda_j e^{S_i}Be^{S_j})_{1 \le i,j \le I},$$ for $S_1,\cdots,S_I \in \mathrm{Sym}(n)$ such that $\sum_{i=1}^I \lambda_i S_i \ge 0.$  Then, $$\det(I_n+B) \le \prod_{i=1}^I \det(I_n+A_i)^{\lambda_i},$$ and there is equality if and only if $A_1=\cdots=A_I=B$ such that $\mathrm{Ran}(B) \subset \cap_{i=1}^I \ker(S_i).$ 
    
    \hfill \\
    More generally, if $A_1,\cdots,A_I,B \in \mathrm{Sym}_{\ge 0}(n)$ such that $\mathrm{diag}(\lambda_1A_1,\cdots,\lambda_I A_I) \ge (\lambda_i \lambda_j e^{S_i}Be^{S_j})_{1 \le i,j \le I}$, where $S_1,\cdots,S_I \in \mathrm{Sym}(n)$ and $$\Gamma \le \mathcal O_\lambda(e^{S_1}\Sigma_1e^{S_1},\cdots,e^{S_I}\Sigma_Ie^{S_I}),$$ then $$\det(I_n+\Gamma B) \le \prod_{i=1}^I \det(I_n+\Sigma_i A_i)^{\lambda_i},$$ and there is equality if and only if \begin{equation} \label{Gaussian_Eq_Cases_L0(B)Thm} A_i=e^{S_i}Be^{S_i} \text{ such that }\mathrm{Ran}(B) \subset \bigcap_{i=1}^I \ker(e^{S_i}\Sigma_i e^{S_i}-\Gamma).\end{equation}
\end{proposition}

\begin{proof}
     Let $E=\ker(B)^\perp$ and $$C=B^{1/2}|_{E}\in \R^{n \times \dim(E)},$$ a full-rank matrix such that $$CC^*=B, C^*C=B|_E>0.$$ Setting $$W=\begin{pmatrix} \lambda_1 e^{S_1}C  \\
     \vdots \\
     \lambda_I e^{S_I}C \end{pmatrix}$$ then $$WW^*=(\lambda_i \lambda_j e^{S_i}Be^{S_j})_{1 \le i,j \le I} \text{ so that } \mathrm{diag}(\lambda_1A_1,\cdots,\lambda_IA_I) \ge WW^*,$$ which implies by the generalized Schur complement, $$W^*\mathrm{diag}\Big(\frac1{\lambda_1}A_1^\dagger, \cdots,\frac1{\lambda_I}A_I^\dagger \Big)W \le I_E,$$ where $A^\dagger$ denotes the Moore-Penrose inverse. Moreover, $\mathrm{diag}(\lambda_1 A_1,\cdots,\lambda_I A_I) \ge WW^*$ implies that $$\ker(\mathrm{diag}(\lambda_1 A_1,\cdots,\lambda_I A_I)) \subset \ker(W^*) \Longrightarrow \mathrm{Ran}(W) \subset \mathrm{Ran}(\mathrm{diag}(\lambda_1 A_1,\cdots,\lambda_I A_I)),$$ which proves that $\mathrm{Ran}(e^{S_i}C) \subset \mathrm{Ran}(A_i)$ for any $1 \le i \le I.$
     \hfill \\
     Thus, we may define $X_i=(e^{S_i}C)^*A_i^\dagger (e^{S_i}C) \in \R^{\dim(E) \times \dim(E)}$ and then, since $e^{S_i}C$ is injective with range contained in $\mathrm{Ran}(A_i)$, $X_i$ is invertible.
     
    \hfill \\
    In addition, as $W^*\mathrm{diag}\Big(\frac1{\lambda_1}A_1^\dagger, \cdots,\frac1{\lambda_I}A_I^\dagger \Big)W \le I_E,$ we have that $$\sum_{i=1}^I \lambda_i X_i \le I_E.$$ Finally, notice that $$A_i \ge (e^{S_i}C)X_i^{-1}(e^{S_i}C)^*.$$ Indeed, $$X_i=\big[(A_i^\dagger)^{1/2}e^{S_i}C \big ]^* \big[ (A_i^\dagger)^{1/2}e^{S_i}C\big] \quad \text{ for any } 1 \le i \le I,$$ and $$\begin{array}{lll}
    \displaystyle A_i^{1/2}\big[(A_i^\dagger)^{1/2}e^{S_i}C\big]=e^{S_i}C \Longrightarrow (e^{S_i}C)X_i^{-1}(e^{S_i}C)^*= \\ \\
    \displaystyle A_i^{1/2}\big[(A_i^\dagger)^{1/2}e^{S_i}C\big]\Big( \big[(A_i^\dagger)^{1/2}e^{S_i}C \big ]^* \big[ (A_i^\dagger)^{1/2}e^{S_i}C\big] \Big)^{-1}\big[(A_i^\dagger)^{1/2}e^{S_i}C\big]^* A_i^{1/2} \le A_i.
    \end{array}$$ 
    Moreover, applying Lemma \ref{Aux_Lemma} with $P_i=S_i,Q_i=-\frac12 \log X_i$ and using here the concavity of the logarithm to get $$\sum_{i=1}^I \lambda_i \log(X_i) \le \log \Big(\sum_{i=1}^I \lambda_i X_i\Big) \le 0,$$ yields that $$\begin{array}{lll}
    \displaystyle \det(I_n+B)=\det(I_E+C^*C) &\le& \displaystyle  \prod_{i=1}^I \det(I_E+(e^{P_i}Ce^{Q_i})^*(e^{P_i}Ce^{Q_i}))^{\lambda_i} \\ \\
    &\le& \displaystyle  \prod_{i=1}^I \det(I_n+(e^{S_i}C)X_i^{-1}(e^{S_i}C)^*)^{\lambda_i} \\ \\
    &\le& \displaystyle \prod_{i=1}^I \det(I_n+A_i)^{\lambda_i}.
    \end{array}$$ If there is equality, we have that $$\sum_{i=1}^I \lambda_i Q_i=0 \Longrightarrow \sum_{i=1}^I \lambda_i \log(X_i)=0,$$ hence by taking traces, $$0=\sum_{i=1}^I \lambda_i \log \det(X_i) \le \log \det \Big(\sum_{i=1}^I \lambda_i X_i\Big) \le 0,$$ by concavity of $\log \det$ and as $\sum_{i=1}^I \lambda_i X_i \le I_E.$ Using the strict concavity of $\log \det$, $$X_1=\cdots=X_I=\sum_{i=1}^I \lambda_i X_i=I_E.$$ In particular, this implies $Q_i=0$ for any $1 \le i \le I$, so that $$P_iC+CQ_i=0 \Longrightarrow S_iC=0 \Longrightarrow \mathrm{Ran}(B) \subset \ker(S_i),$$ for any $1 \le i \le I.$ Finally, $A_i=(e^{S_i}C)X_i^{-1}(e^{S_i}C)^*=CC^*=B$, since $e^{S_i}C=C$ for any $1 \le i \le I.$ We conclude, as desired, that $$A_1=\cdots=A_I=B, \quad \mathrm{Ran}(B) \subset \bigcap_{i=1}^I \ker(S_i).$$
    \hfill \\
    For the more general version, defining $$R_i=e^{S_i}\Sigma_i e^{S_i}, \quad M=\mathcal O_\lambda(R_1,\cdots,R_I),$$ we set $$T_i=\frac12 \log (M^{-1/2}R_iM^{-1/2}) \quad \forall 1 \le i \le I.$$ Then, $\sum_{i=1}^I \lambda_i T_i=0$ and writing $L_i=e^{T_i}M^{1/2}e^{-S_i}$, $$L_ie^{S_i}=e^{T_i}M^{1/2} \text{ and } L_i^TL_i=e^{-S_i}M^{1/2}e^{2T_i}M^{1/2}e^{-S_i}=\Sigma_i,$$ for any $1 \le i \le I.$ If we conjugate $$\mathrm{diag}(\lambda_1 A_1,\cdots,\lambda_I A_I) \ge (\lambda_i \lambda_j e^{S_i}Be^{S_j})_{1 \le i,j \le I}$$ by $\mathrm{diag}(L_1,\cdots,L_I)$, one has $$\mathrm{diag}(\lambda_1 L_1A_1L_1^T,\cdots,\lambda_I L_IA_IL_I^T) \ge (\lambda_i \lambda_j e^{T_i}M^{1/2}BM^{1/2}e^{T_j})_{i,j}.$$ Hence, $$\det(I_n+M^{1/2}BM^{1/2}) \le \prod_{i=1}^I \det(I_n+L_iA_iL_i^T)^{\lambda_i}.$$ As $B^{1/2}\Gamma B^{1/2} \le B^{1/2}MB^{1/2}$, $$\det(I_n+B\Gamma)=\det(I_n+B^{1/2}\Gamma B^{1/2}) \le \det(I+BM) \le \prod_{i=1}^I \det(I_n+\Sigma_i A_i)^{\lambda_i}.$$ 
    \hfill \\
    For the equality cases, by what was previously shown, we need $$M^{1/2}BM^{1/2}=L_iA_iL_i^T \text{ and } T_i M^{1/2}BM^{1/2}=0 \quad \forall 1 \le i \le I.$$ Then, $$M^{-1/2}(R_i-M)BM^{1/2}=(e^{2T_i}-I_n)M^{1/2}BM^{1/2}=0,$$ which implies that $(R_i-M)B=0$ for any $1 \le i \le I.$ 
    In addition, $$e^{-T_i}L_iA_iL_i^Te^{-T_i}=M^{1/2}BM^{1/2} \Longrightarrow e^{-S_i}A_ie^{-S_i}=B,$$ for any $1 \le i \le I.$ 
    Finally, $\det(I_n+\Gamma B)=\det(I_n+M^{1/2}BM^{1/2})$ is equivalent to $B^{1/2}(M-\Gamma)B^{1/2}=0$ or in other words, $(M-\Gamma)B=0.$ 
    
    \hfill \\
    For the sufficiency of \cref{Gaussian_Eq_Cases_L0(B)Thm}, it suffices to observe that $$\det(I_n+\Sigma_i A_i)=\det(I_n+R_iB)=\det(I_n+\Gamma B).$$ 
\end{proof}

\hfill \\
This allows us to characterize the Gaussian equality cases of $$\prod_{i=1}^I f_i^{\lambda_i}(x_i) \le h \Big(\sum_{i=1}^I \lambda_i e^{S_i}x_i\Big) \Longrightarrow \prod_{i=1}^I \Big(\int_{\R^n}f_i\, d\gamma_{\Sigma_i} \Big)^{\lambda_i} \le \int_{\R^n}h\, d\gamma_\Gamma.$$ We now seek to understand the equality cases of \cref{Multiple_BThm}, whose proof will follow the lines of the derivation of equality cases for the strong Gaussian correlation in \cite{Strong_Gaussian_Correlation}.

\hfill \\
For the rest of this subsection, fix an extremizing tuple $(f_i)_{1 \le i \le I}$ such that $f_i$ is log-concave and centered with respect to $\gamma_{\Sigma_i}$ and $h$ such that it is log-concave. We may normalize $f_1,\cdots,f_I$ and $h$ so that $$\int_{\mathbb{R}^n}f_i\, d\gamma_{\Sigma_i}=\int_{\mathbb{R}^n}h\, d\gamma_\Gamma=1 \quad \forall 1 \le i \le I.$$ 
Recall the definition of rescaled self-convolution $$\mathrm{Conv}_{\gamma_\Sigma}(f)(x)=\int_{\mathbb{R}^n}f\Big(\frac{x+y}{\sqrt2}\Big) f\Big(\frac{x-y}{\sqrt 2}\Big)\, d\gamma_\Sigma(y).$$

\begin{lemma}
For any $x_1,\cdots,x_I \in \mathbb{R}^n$, $$\prod_{i=1}^I \mathrm{Conv}_{\gamma_{\Sigma_i}}(f_i)^{\lambda_i}(x_i) \le \mathrm{Conv}_{\gamma_\Gamma}(h)\Big(\sum_{i=1}^I \lambda_i e^{S_i}x_i\Big).$$ 
\end{lemma}

\begin{proof}
    Fixing $x_1,\cdots,x_I \in \mathbb{R}^n,$ we set $$z=\sum_{i=1}^I \lambda_i e^{S_i}x_i,$$ we may define the functions $$F_{i,x_i}(u_i)=f_i\Big( \frac{x_i+u_i}{\sqrt 2}\Big)f_i\Big(\frac{x_i-u_i}{\sqrt 2}\Big) \quad \forall 1 \le i \le I,$$ and $$H_z(w)=h\Big(\frac{z+w}{\sqrt 2}\Big)h\Big(\frac{z-w}{\sqrt 2}\Big).$$ It is clear that $F_{i,x_i}$ and $H_z$ are even and log-concave for any $1 \le i \le I$. The pointwise inequality involving $f_1,\cdots,f_I$ and $h$ is easily seen to be transferred to the same pointwise inequality for $F_{1,x_1},\cdots,F_{I,x_I}$ and $H_z.$ Applying \cref{FRBL_MultipleBthm} proves the desired claim. 
\end{proof}
\hfill \\
By the local central limit theorem, $$\mathrm{Conv}_{\gamma_{\Sigma_i}}^N(f_i) \rightarrow_{N \rightarrow +\infty} \frac{\gamma_{A_i}}{\gamma_{\Sigma_i}}, \text{ where } A_i=\mathrm{Cov}(f_i \gamma_{\Sigma_i})>0,$$ in pointwise convergence.

\begin{lemma}
    Under equality of \cref{FRBL_MultipleBthm}, one has $$\mathrm{bar}_{\gamma_\Gamma}(h)=0.$$
\end{lemma}

\begin{proof}
    Assume by contradiction that $m=\mathrm{bar}_{\gamma_\Gamma}(h) \neq 0.$ Then, the probability measure $\mathrm{Conv}_{\gamma_\Gamma}^N(h)\, d\gamma_\Gamma$ has mean $2^{N/2} m$ and a covariance matrix, that is independent of $N.$ Applying iteratively the relation of the previous lemma, we get that $$\mathrm{Conv}_{\gamma_\Gamma}^N(h)\Big(\sum_{i=1}^I \lambda_i e^{S_i}x_i \Big) \ge \prod_{i=1}^I \mathrm{Conv}_{\gamma_{\Sigma_i}}^N(f_i)(x_i)^{\lambda_i},$$ or in other words, considering $x_i=e^{-S_i}w$, $$\mathrm{Conv}_{\gamma_\Gamma}^N(h)(w) \ge \prod_{i=1}^I \mathrm{Conv}_{\gamma_{\Sigma_i}}^N(f_i)(e^{-S_i}w)^{\lambda_i}.$$ The limits of $\mathrm{Conv}_{\gamma_{\Sigma_i}}^N(f_i)$ are strictly positive, so that, defining $$u=\frac m{|m|}, H_-=\{ w: \langle u,w\rangle \le 0\},$$ one has, on the one hand, that $$\liminf_{N \rightarrow +\infty} \int_{H_-}\mathrm{Conv}_{\gamma_\Gamma}^N(h) \, d\gamma_\Gamma \ge \int_{H_-}\prod_{i=1}^I \big[(\gamma_{A_i}/\gamma_{\Sigma_i})(e^{-S_i}w)\big]^{\lambda_i}\, d\gamma_\Gamma(w)>0.$$ On the other hand, by Chebyshev's inequality, $$\int_{H_-}\mathrm{Conv}_{\gamma_\Gamma}^N(h)\, d\gamma_\Gamma \le \frac{\langle \mathrm{Cov}(h \gamma_\Gamma)u,u\rangle}{2^N |m|^2} \rightarrow_{N \rightarrow +\infty} 0.$$ This yields the desired contradiction. 
\end{proof}

\hfill \\
We may then also apply the local central limit theorem to $\mathrm{Conv}_{\gamma_\Gamma}^N(h)$, that converges pointwise to $$\frac{\gamma_{\Sigma_h}}{\gamma_\Gamma}, \quad \Sigma_h=\mathrm{Cov}(h \gamma_\Gamma)>0.$$

\begin{corollary}
    Let $((f_i)_{1 \le i \le I},h)$ be a normalized extremizing tuple for \cref{FRBL_MultipleBthm}. Then, $h$ has Gaussian barycenter $0$ and the centered Gaussian relative densities with covariances $$A_i=\mathrm{Cov}(f_i\gamma_{\Sigma_i}), \quad \Sigma_h=\mathrm{Cov}(h \gamma_\Gamma),$$ satisfy $$\prod_{i=1}^I \Big( \frac{\gamma_{A_i}}{\gamma_{\Sigma_i}}\Big)^{\lambda_i}(x_i) \le \frac{\gamma_{\Sigma_h}}{\gamma_\Gamma}\Big(\sum_{i=1}^I \lambda_i e^{S_i}x_i\Big) \quad \forall x_1,\cdots,x_I \in \mathbb{R}^n.$$ 
\end{corollary}
\hfill \\
Note that if $f$ is log-concave, the probability measure $$\mu_{f,\Sigma}=\frac{f\, d\gamma_\Sigma}{\int f\, d\gamma_\Sigma} \text{ is more log-concave than } \gamma_\Sigma.$$ Applying the Brascamp-Lieb variance inequality implies that $\mathrm{Cov}(\mu_{f,\Sigma}) \le \Sigma.$ Thus, setting $$B_i=A_i^{-1}-\Sigma_i^{-1} \ge 0, \quad B=\Sigma_h^{-1}-\Gamma^{-1} \ge 0,$$ and the densities in the previous corollary are $$\frac{d\gamma_{A_i}}{d\gamma_{\Sigma_i}}(x)=\det(I_n+\Sigma_i B_i)^{1/2}g_{B_i}(x), \quad \frac{d\gamma_{\Sigma_h}}{d\gamma_\Gamma}(x)=\det(I_n+\Gamma B)^{1/2}g_B(x).$$ Evaluating the limiting pointwise inequality at the origin yields $$\det(I_n+\Gamma B)^{1/2} \ge \prod_{i=1}^I \det(I_n+\Sigma_i B_i)^{\lambda_i/2},$$ while \cref{Gaussian_Case_MultipleBThm} gives the reverse inequality.
By the characterization of equality cases for Gaussian functions, we conclude that \begin{equation} \label{Forced_Equalities_Covariances} B_i=e^{S_i}Be^{S_i} \quad \forall 1 \le i \le I, \quad \mathrm{Ran}(B) \subset \bigcap_{i=1}^m \ker(e^{S_i}\Sigma_i e^{S_i}-\Gamma).\end{equation}

\hfill \\
The following splitting lemma turns this matrix information into information about the original functions. It is a consequence of \cite{Shohei_Hiroshi2} and of a linear change of variables. 

\begin{lemma} \label{Splitting_Lemma}
    Let $f:\R^n \rightarrow [0;+\infty)$ be log-concave, $$\int f\, d\gamma_\Sigma=1, \int xf(x)\, d\gamma_\Sigma(x)=0.$$ Then, there exists a log-concave function $\psi: \big[\ker \big(\mathrm{Cov}(f\gamma_\Sigma)^{-1}-\Sigma^{-1} \big) \big]^\perp \rightarrow [0;+\infty)$, such that $$f(x)=\psi \Big(P_{\big[\ker \big(\mathrm{Cov}(f\gamma_\Sigma)^{-1}-\Sigma^{-1} \big) \big]^\perp}x\Big) \text{ for a.e. } x \in \R^n.$$
\end{lemma}

\hfill \\
Applying the previous lemma and setting $F=\mathrm{Ran}(B)$, $$h(x)=\psi_0(P_Fx), \quad f_i(x_i)=\psi_i(P_Fe^{S_i}x_i),$$ for some log-concave $\psi_i: F \rightarrow [0;+\infty), \; 0 \le i \le I.$ We used here that $$\ker(B_i)=e^{-S_i}\ker(B)=e^{-S_i}F^\perp.$$ 
By \cref{Forced_Equalities_Covariances}, $$P_F\big( e^{S_i}\Sigma_ie^{S_i}\big)P_F\big|_F=\Gamma_F \text{ for any } 1 \le i \le I.$$ Hence, the pointwise condition becomes $$\psi_0\Big(\sum_{i=1}^I \lambda_i y_i\Big) \ge \prod_{i=1}^I \psi_i(y_i)^{\lambda_i},$$ while the integral equality is \begin{equation} \label{Equality_ArbitraryCov_PL} \int_F \psi_0 \, d\gamma_{\Gamma_F}=\prod_{i=1}^I \Big(\int_F \psi_i\, d\gamma_{\Gamma_F}\Big)^{\lambda_i}.\end{equation}

\hfill \\
Therefore, the original equality problem reduces exactly to equality in the Gaussian Pr\'ekopa-Leindler on $F.$ For the equality cases of Gaussian Pr\'ekopa-Leindler, we need to use the characterization of equality cases of \cite[Thm 4.1]{Kolesnikov_Werner} for the Pr\'ekopa-Leindler with multiple functions. 

\begin{lemma} \label{Gaussian_PL_EqCases}
    Let $f_1,\cdots,f_I, h: \R^n \rightarrow [0;+\infty)$ be log-concave functions such that $$0<\int_{\mathbb{R}^n}f_i\, d\gamma_n, \; \int_{\mathbb{R}^n}h\, d\gamma_n<+\infty.$$ Assume that $$h\Big(\sum_{i=1}^I \lambda_i x_i\Big) \ge \prod_{i=1}^I f_i(x_i)^{\lambda_i} \quad \forall x_1,\cdots,x_I \in \R^n,$$ and $$\int_{\R^n} h\, d\gamma_n=\prod_{i=1}^I \Big(\int_{\R^n}f_i \, d\gamma_n \Big)^{\lambda_i}.$$ Then, there exist $a_1,\cdots,a_I>0$ and a non-zero log-concave $\phi$, such that $$f_i=a_i \phi \quad \forall 1 \le i \le I \text{ and } h=\prod_{i=1}^I a_i^{\lambda_i} \phi \text{ almost everywhere.}$$ If one of the $f_i's$ is centered w.r.t $\gamma$, so is $\phi.$
\end{lemma}

\begin{proof}
    It suffices to observe that, for any $x_1,\cdots,x_I \in \R^n$, $$\begin{array}{lll}
    \displaystyle \prod_{i=1}^I (f_i \gamma_n)(x_i)^{\lambda_i}=(2\pi)^{-n/2}e^{-\frac12 \sum_{i=1}^I \lambda_i |x_i|^2} \prod_{i=1}^I f_i(x_i)^{\lambda_i} \le (2\pi)^{-n/2}e^{-\frac12\Big|\sum_{i=1}^I \lambda_i x_i\Big|^2} \\ \\ \displaystyle e^{-\frac12 \Big(\sum_{i=1}^I \lambda_i |x_i|^2-\Big|\sum_{i=1}^I \lambda_i x_i\Big|^2\Big)}h\Big(\sum_{i=1}^I \lambda_i x_i\Big)=e^{-\frac14 \sum_{i,j=1}^I \lambda_i \lambda_j |x_i-x_j|^2}(h\gamma_n)\Big(\sum_{i=1}^I \lambda_i x_i\Big) \\ \\ 
    \displaystyle \le (h\gamma_n)\Big(\sum_{i=1}^I \lambda_i x_i\Big).
    \end{array}$$
    Hence, by \cite[Thm 4.1]{Kolesnikov_Werner}, $$(f_i\gamma_n)(x_i)=a_i F(x_i-v_i) \quad \forall 1 \le i \le I, \quad (h\gamma_n)(x)=\prod_{i=1}^I a_i^{\lambda_i}F\Big(x-\sum_{i=1}^I \lambda_i v_i\Big) \text{ a.e.},$$ where $F$ is log-concave and $a_i>0, v_i \in \R^n$ for any $1 \le i \le I.$ But $$\begin{array}{lll}
    \ds F(x)=\frac{\prod_{i=1}^I (f_i\gamma_n)(x+v_i)^{\lambda_i}}{\prod_{i=1}^I a_i^{\lambda_i}} &\le& \ds e^{-\frac14\sum_{i,j=1}^I \lambda_i \lambda_j |v_i-v_j|^2}\frac{(h\gamma_n)(x+\sum_{i=1}^I \lambda_i v_i)}{\prod_{i=1}^I a_i^{\lambda_i}} \\ 
    \nm 
    &=& \ds e^{-\frac14\sum_{i,j=1}^I \lambda_i \lambda_j |v_i-v_j|^2} F(x).
    \end{array}$$
    This forces $v_i=v_j$ for any $i,j$ and thus, setting $\phi=\frac{F(\cdot-v_1)}{\gamma_n}$ concludes the proof of the lemma. 
\end{proof}
\hfill \\
By change of variables, this provides the same equality cases for \cref{Equality_ArbitraryCov_PL}. Combining the previous lemmas and checking that the necessity of the special form of functions attaining equality is in fact sufficient yields the following theorem.

\begin{theorem} \label{Eq_Cases_MultipleBThm}
Set $I \ge 2, \lambda \in \Delta_I^\circ$ and $S_1,\cdots,S_I \in \mathrm{Sym}(n), \Gamma, \Sigma_1,\cdots,\Sigma_I \in \mathrm{Sym}_{>0}(n)$ such that $$\Gamma \le \mathcal O_\lambda(e^{S_1}\Sigma_1e^{S_1},\cdots,e^{S_I}\Sigma_I e^{S_I}).$$ Define $$E=\bigcap_{i=1}^I \ker(e^{S_i}\Sigma_ie^{S_i}-\Gamma).$$
    Let $f_1,\cdots,f_I,h:\R^n \rightarrow [0;+\infty)$ be log-concave functions with $$0<\int_{\mathbb{R}^n}f_i\, d\gamma_{\Sigma_i}, \; \int_{\mathbb{R}^n}h\, d\gamma_\Gamma<+\infty.$$ Assume that $f_i$ are centered w.r.t $\gamma_{\Sigma_i}$ and $$h\Big(\sum_{i=1}^I \lambda_i e^{S_i}x_i\Big) \ge \prod_{i=1}^I f_i(x_i)^{\lambda_i}, \quad \forall x_1,\cdots,x_I \in \R^n.$$ Then, $$\int_{\R^n}h \, d\gamma_\Gamma \ge \prod_{i=1}^I \Big(\int_{\R^n}f_i\, d\gamma_{\Sigma_i} \Big)^{\lambda_i},$$ and there is equality iff there exist constants $a_1,\cdots,a_I>0$ and a non-zero log-concave function $$\phi: E \longrightarrow [0;+\infty), \text{ centered with respect to } \gamma_{\Gamma_E},$$ so that $$f_i(x)=a_i \phi(P_Ee^{S_i}x) \quad \forall 1 \le i \le I, h(x)=\prod_{i=1}^I a_i^{\lambda_i} \phi(P_Ex),$$ for a.e. $x \in \R^n.$ 
\end{theorem}

\hfill \\
Restricting the theorem to characteristic functions provides the equality cases of the corresponding geometric inequality.  

\begin{corollary}
In the same setting as the above theorem, assume that $K_1,\cdots,K_I \subset \R^n$ are convex sets, with non-empty interior and $$\int_{K_i}x\, d\gamma_{\Sigma_i}(x)=0 \quad \forall 1 \le i \le I.$$ Then, $$\gamma_\Gamma\Big(\sum_{i=1}^I \lambda_i e^{S_i}K_i\Big) \ge \prod_{i=1}^I \gamma_{\Sigma_i}(K_i)^{\lambda_i}.$$ Equality holds iff there exists a centered (w.r.t $\gamma_{\Gamma_E}$) convex set $C \subset E$, such that $$e^{S_i}K_i=C+E^\perp \quad \forall 1 \le i \le I,$$ up to null-sets. 
\end{corollary}

\hfill \\
\cref{Multiple_BThm} naturally goes in the direction of the (B)-theorem (see \cite{Cordero_Fradelizi_Maurey,Cordero_Rotem,Cordero_Eskenazis}). More recently, \cite{Faifman_Polt} asked whether the functions $$t \mapsto \gamma_n(e^tK) \text{ or more generally, } t \mapsto \gamma_n(e^{tS}K) \text{ for a symmetric matrix } S$$ are log-concave for $K \subset \R^n$ a convex body whose Lebesgue barycenter is at the origin. 

\hfill \\
In fact, we answer this question negatively by providing a counterexample already in the three-dimensional case. Let us note that the one-dimensional case is trivially true by the usual Gaussian (B)-theorem of \cite{Cordero_Fradelizi_Maurey} as centering and symmetry coincide on the real line. However, it is not clear whether the question has a positive or negative answer in the planar case. 

\begin{theorem} \label{Counterexample_FailureLogCcv}
    There exists a convex body $K \subset \mathbb{R}^3$ centered with respect to the Lebesgue measure, such that $t \mapsto \gamma_3(e^tK)$ is not log-concave. 
\end{theorem}

\begin{remark}
    More generally, for $n>3$, the convex bodies $K\times[-L,L]^{n-3}$ have vanishing Lebesgue barycenter and provide counterexamples for all sufficiently large $L$, since
$$
\frac{d^2}{dt^2}\log\gamma_1([-Le^t,Le^t])\big|_{t=0}
\rightarrow 0
\qquad\text{as }L\rightarrow +\infty.$$
\end{remark}

\begin{proof}[Proof of \cref{Counterexample_FailureLogCcv}]
Define $$K_R=\Big \{(x,y) \in \R \times \R^2: -R \le x \le 3R, \quad |y| \le \frac{3R-x}{R^3} \Big\}, \quad F_R(t)=\log \gamma_3(e^tK_R).$$ 
First, the Lebesgue barycenter of $K_R$ is $0$, as $$\int_{K_R}y\, dx\, dy=0 \text{ by symmetry, }$$ while $$\begin{array}{lll}
\ds \int_{K_R}x\, dx dy &=& \ds \frac \pi{R^6}\int_{-R}^{3R} x(3R-x)^2\, dx \\ 
\nm
&=& \ds \frac{\pi}{R^6}\int_0^{4R}x^2(3R-x)\, dx=\frac{\pi}{R^6}(64R^4-\frac14 (4R)^4)=0.
\end{array}$$ The goal is to show that $$F_R''(0)>0 \text{ for every sufficiently large } R.$$ Let $$\mu_R=\frac{\mathbbm{1}_{K_R}(z)}{\gamma_3(K_R)}d\gamma_3, \quad Z=(X,Y) \sim \mu_R,$$ and $$Q=|Z|^2=X^2+|Y|^2.$$ Then, it is standard (see \cite{Cordero_Fradelizi_Maurey}) that $$F_R''(0)=\mathbb{E}_{\mu_R}(Q^2)-(\mathbb{E}_{\mu_R} Q)^2-2\mathbb{E}_{\mu_R}Q.$$ We therefore need to compute the first two moments of $Q$ sufficiently accurately. 
\hfill \\
First, consider for any $L>0$, the probability measure $$d\hat{\nu}_L=\frac{(L-x)^2}{L^2+1}\frac1{\sqrt{2\pi}}e^{-\frac{|x|^2}{2}}\, dx,$$ and denote $\varphi(x)=\frac1{\sqrt{2\pi}}e^{-\frac{|x|^2}{2}}.$ Using the standard Gaussian moments $$\int_{\R} x^2\varphi(x)\, dx=1, \quad \int_{\R} x^4 \varphi(x)\, dx=3, \quad \int_{\R}x^6 \varphi(x)\, dx=15,$$ we have that $$\int_{\R} x^2(L-x)^2\varphi(x)\, dx=L^2+3,$$ and $$\int x^4(L-x)^2 \varphi(x)\, dx=3L^2+15,$$ so that $$\hat m_2(L)=\mathbb{E}_{\hat \nu_L}(X^2)=\frac{L^2+3}{L^2+1}, \quad \hat m_4(L)=\mathbb{E}_{\hat \nu_L}(X^4)=\frac{3L^2+15}{L^2+1}.$$ Then, $$\mathrm{Var}_{\hat \nu_L}(X^2)-2\mathbb{E}_{\hat \nu_L}(X^2)=\frac{4L^2}{(L^2+1)^2}.$$ 
The remaining work is to justify replacing the Gaussian moments on $K_R$ by the ones just computed with $L=3R.$ Set $$I_R=[-R;3R], \quad L=3R,$$ and introduce the truncated probability measure $$d\nu_R(x)=\frac{(L-x)^2\varphi(x)\mathbbm{1}_{I_R}(x)}{N_R}\, dx, \text{ where } N_R=\int_{I_R}(L-x)^2 \varphi(x)\, dx.$$ Let us now compare $\nu_R$ with $\hat{\nu}_L.$ On $\R\setminus I_R$, $|x| \ge R$ so that $(L-x)^2 \le (3R+|x|)^2 \le 16x^2$, which implies that $$\int_{\R \setminus I_R}x^j(L-x)^2 \varphi(x)\, dx \le 32\int_R^{+\infty}x^{j+2}\varphi(x)\, dx=32\Big(R^{j+1}\varphi(R)+(j+1)\int_R^{+\infty}x^j \varphi(x)\, dx \Big),$$ so that $$\int_{\R\setminus I_R}x^j(L-x)^2 \varphi(x)\, dx=O(R^{j+1}\varphi(R)) \text{ for } j=0,2,4,$$ also using the standard estimate $$\int_R^{+\infty}\varphi(x)\, dx \le \frac{\varphi(R)}{R}.$$ Therefore, $$\mathbb{E}_{\nu_R}(X^2)=\hat{m}_2(L)+O(R\varphi(R)), \quad \mathbb{E}_{\nu_R}(X^4)=\hat m_4(L)+O(R^3\varphi(R)).$$
\hfill \\
Set $\eta_R$ the distribution of $X$ when $(X,Y) \sim \mu_R$. We want to compare $\eta_R$ with $\nu_R.$ By polar coordinates, the density of $\eta_R$ is proportional to $$\varphi(x)\Big(1-e^{-(3R-x)^2/2R^6}\Big)\mathbbm{1}_{I_R}(x).$$ Note that $$1-e^{-(3R-x)^2/2R^6}=\frac{(3R-x)^2}{2R^6}\int_0^1 e^{-s(3R-x)^2/2R^6}\, ds,$$ and as $-R \le x \le 3R$, $$e^{-8/R^4}\le \int_0^1 e^{-s(3R-x)^2/2R^6}\, ds \le 1.$$ Hence, $$e^{-8/R^4}\le \frac{d\eta_R}{d\nu_R}(x)=\frac{\int_0^1 e^{-s(3R-x)^2/2R^6}\, ds}{\mathbb{E}_{\nu_R}(\int_0^1 e^{-s(3R-\cdot)^2/2R^6}\, ds)} \le e^{8/R^4}.$$ This means that for every non-negative $g$, $$\Big| \mathbb{E}_{\eta_R}(g)-\mathbb{E}_{\nu_R}(g)\Big| \le (e^{8/R^4}-1)\mathbb{E}_{\nu_R}(g),$$ which means in particular, $$\mathbb{E}_{\mu_R}(X^2)=\hat{m}_2(L)+O(R^{-4})+O(R\varphi(R)) \Big(=O(1)\Big),$$ and $$\mathbb{E}_{\mu_R}(X^4)=\hat m_4(L)+O(R^{-4})+O(R^3\varphi(R)).$$ 
\hfill \\
Recall that $Q=X^2+|Y|^2$, with $0 \le |Y|^2 \le 16R^{-4}$ so that $$0 \le \mathbb{E}_{\mu_R}(Q)-\mathbb{E}_{\mu_R}(X^2) \le 16R^{-4},$$ and $$\begin{array}{lll}
0 &\le& \ds \mathbb{E}_{\mu_R}(Q^2)-\mathbb{E}_{\mu_R}(X^4)=2\mathbb{E}_{\mu_R}(X^2|Y|^2)+\mathbb{E}_{\mu_R}(|Y|^4) \\ \\
&\le& \ds \frac{32}{R^4}\mathbb{E}_{\mu_R}(X^2)+\frac{256}{R^8}=O(R^{-4}).
\end{array}$$ As $R\varphi(R)$ decays faster than any power of $R$, we conclude that $$\mathbb{E}_{\mu_R}(Q)=\hat m_2(L)+O(R^{-4}), \quad \mathbb{E}_{\mu_R}(Q^2)=\hat m_4(L)+O(R^{-4}).$$ Therefore, $$F_R''(0)=\frac{36R^2}{(9R^2+1)^2}+O(R^{-4})=\frac{4}{9R^2}+O(R^{-4})>0,$$ for $R$ sufficiently large. 
\end{proof}

\hfill \\
\begin{remark}
    The claim fails even when the assumption of the convex set $K$ having zero Lebesgue barycenter is replaced by the (possibly more natural) assumption of zero Gaussian barycenter. Fix a Lebesgue-centered convex body $K$ for which $$F_K''(0)>0, \quad F_K(t)=\log \gamma_n(e^tK).$$ For $\epsilon>0$ small enough, define $$G(\epsilon,a)=\int_K (x+a)e^{-\epsilon^2|x+a|^2/2}\, dx.$$ Then, $$G(0,0)=0, \quad D_aG(0,0)=|K|I_n,$$ so that by the implicit function theorem, there exists $a_\epsilon \rightarrow 0$ with $G(\epsilon,a_\epsilon)=0.$ Hence, $$\tilde K_\epsilon=\epsilon(K+a_\epsilon)$$ has Gaussian barycenter zero. Then, $$F_{\tilde K_\epsilon}''(\log(1/\epsilon))=F_{K+a_\epsilon}''(0)>0,$$ by continuity of $a \mapsto F_{K+a}''(0)$ and as at $a=0$, it is positive by assumption. Hence, global (B) log-concavity also may fail for a Gaussian centered convex body. 
\end{remark}

\hfill \\
On the positive side, if we take account of all translations, we indeed recover log-concavity. In other words, optimizing over translations restores the (B)-property, while fixed centering can fail in dimensions at least three. 

\begin{corollary} \label{Translated_BThm}
    For every convex set $K \subset \R^n$ and symmetric matrix $S \in \mathrm{Sym}(n)$, $$t \mapsto \sup_{a \in \R^n}\gamma_n(e^{tS}(K+a)) \text{ is log-concave.}$$
\end{corollary}

\begin{remark}
    If $K$ is supposed to be symmetric, this recovers by Anderson's theorem (see \cite{Anderson,Gardner_BM}) the (B)-theorem for symmetric sets of \cite{Cordero_Fradelizi_Maurey}. 
\end{remark}

\begin{proof}
Let us first assume that $K$ is a convex body, \textit{i.e.} a compact convex subset of $\R^n$ with non-empty interior. 
It is then clear that for any $t \in \R$, $$F_K(t):=\sup_{a \in \mathbb{R}^n}\gamma_n(e^{tS}(K+a)) \text{ is a maximum,}$$ since $\gamma_n(e^{tS}(K+a)) \rightarrow_{|a| \rightarrow +\infty}0$ by boundedness of $K$ and by continuity of the function at $a$, for fixed $t.$

\hfill \\ Let $t_1,t_2 \in \R$ and $\lambda \in (0;1).$ 
We may thus choose maximizing translations $a_1$ at $t_1$ and $a_2$ at $t_2$, \textit{i.e.} $$F_K(t_1)=\gamma_n(e^{t_1S}(K+a_1)), \quad F_K(t_2)=\gamma_n(e^{t_2S}(K+a_2)).$$ Differentiating with respect to translation shows that $$\int_{e^{t_iS}(K+a_i)}x\, d\gamma_n(x)=0 \text{ for } i=1,2.$$ Then, by \cref{Multiple_BThm}, as $$\lambda \big[(1-\lambda)(t_2-t_1)S\big]+(1-\lambda)\big[\lambda(t_1-t_2)S\big]=0,$$ we have that $$\begin{array}{lll}
\ds \gamma_n\Big(e^{(\lambda t_1+(1-\lambda)t_2)S}(K+a_{\lambda})\Big) &=& \ds \gamma_n\Big( \lambda e^{(1-\lambda)(t_2-t_1)S}[e^{t_1S}(K+a_1)] \\ 
\nm
&& \ds +(1-\lambda)e^{\lambda(t_1-t_2)S}[e^{t_2S}(K+a_2)]\Big) \\
\nm 
&\ge& \ds \gamma_n(e^{t_1S}(K+a_1))^\lambda \cdot \gamma_n(e^{t_2S}(K+a_2))^{1-\lambda},
\end{array}$$ where $$a_\lambda=\lambda a_1+(1-\lambda)a_2.$$ This implies that $$F_K(\lambda t_1+(1-\lambda)t_2) \ge F_K(t_1)^\lambda F_K(t_2)^{1-\lambda},$$ which proves the desired log-concavity. 
\hfill \\
For closed convex sets $K$ with non-empty interior, it suffices to consider $$K_m=K \cap mB_2^n \text{ convex bodies for } m \ge 1,$$ large enough. Then, $F_{K_m}$ is log-concave and as $$\gamma_n(e^{tS}(K_m+a)) \uparrow_{m \rightarrow +\infty} \gamma_n(e^{tS}(K+a)),$$ we conclude that $$F_K(t)=\sup_m F_{K_m}(t)=\lim_{m \rightarrow +\infty}F_{K_m}(t),$$ and is hence also log-concave. 
\hfill \\
If $K$ has non-empty interior but is not closed, we may replace it by its closure $\overline K.$ As $\overline K \setminus K \subset \partial \overline K$ and the boundary of a full-dimensional convex set has Lebesgue measure zero, $$\gamma_n(e^{tS}(K+a))=\gamma_n(e^{tS}(\overline K+a)) \quad \forall t \in \R, a \in \R^n \Longrightarrow F_K=F_{\overline K}.$$ Finally, if $K$ has empty interior, it lies in a proper affine subspace and so, $F_K=0.$ 
\end{proof}

\hfill \\

\subsection{Gaussian Borell-Brascamp-Lieb interpolation}

\hfill \\
We now want to consider inequalities of the form $$\gamma_\Gamma\Big(\sum_{i=1}^I \lambda_i e^{S_i}K_i\Big)^{1/n} \ge \mathcal M_p^\lambda(\gamma_{\Sigma_1}(K_1)^{1/n},\cdots,\gamma_{\Sigma_I}(K_I)^{1/n}),$$ for $p \in (0,1]$ and therefore conclude the proof of \cref{Main_Thm_IneqPart}.

\hfill \\
Note that for any $\tilde \lambda \in \Delta_I^\circ$, we have: $$\begin{array}{lll}
    \displaystyle \gamma_\Gamma\Big(\sum_{i=1}^I \lambda_i e^{S_i}K_i\Big)^{1/n} &\ge& \displaystyle \mathcal M_p^\lambda \Big(\gamma_{\Sigma_1}(K_1)^{1/n},\cdots,\gamma_{\Sigma_I}(K_I)^{1/n}\Big)=\Big(\sum_{i=1}^I \tilde \lambda_i \frac{\lambda_i}{\tilde \lambda_i} \gamma_{\Sigma_i}(K_i)^{p/n} \Big)^{1/p}  \\ \\
     &\ge& \displaystyle \Big(\prod_{i=1}^I \Big(\frac{\lambda_i}{\tilde \lambda_i} \gamma_{\Sigma_i}(K_i)^{p/n}\Big)^{\tilde \lambda_i} \Big)^{1/p}=\prod_{i=1}^I \Big(\frac{\lambda_i}{\tilde \lambda_i}\Big)^{\tilde \lambda_i/p} \prod_{i=1}^I \gamma_{\Sigma_i}(K_i)^{\tilde \lambda_i/n},
\end{array}$$ so that $$\gamma_\Gamma \Big(\sum_{i=1}^I \lambda_i e^{S_i}K_i\Big) \ge \prod_{i=1}^I \Big(\frac{\lambda_i}{\tilde \lambda_i}\Big)^{\tilde \lambda_i n/p} \prod_{i=1}^I \gamma_{\Sigma_i}(K_i)^{\tilde \lambda_i},$$ for any $K_1,\cdots,K_I \subset \R^n$ convex with $K_i$ centered with respect to $\gamma_{\Sigma_i}$ for any $1 \le i \le I.$ 

\hfill \\
It then looks natural to consider the corresponding Forward-Reverse Brascamp-Lieb inequality.

\begin{proposition} \label{Computation_towardsBBL}
    Set $p \in (0;1]$ and $\lambda, \tilde \lambda \in \Delta_I^\circ.$ Further, consider $S_1,\cdots,S_I \in \mathrm{Sym}(n)$ and $\Gamma, \Sigma_1,\cdots,\Sigma_I \in \mathrm{Sym}_{>0}(n)$ such that $$\Gamma \le \prod_{i=1}^I \Big(\frac{\tilde \lambda_i}{\lambda_i}\Big)^{2(\frac1p-1)\tilde \lambda_i}\mathcal O_{\tilde \lambda}(e^{S_1}\Sigma_1e^{S_1},\cdots,e^{S_I}\Sigma_I e^{S_I}).$$ If $$h\Big( \sum_{i=1}^I \lambda_i e^{S_i}x_i \Big) \ge \prod_{i=1}^I f_i(x_i)^{\tilde \lambda_i},$$ for $h$ log-concave and $f_i$ log-concave, centered with respect to $\gamma_{\Sigma_i}, 1 \le i \le I$, then $$ \int h\, d\gamma_\Gamma \ge \prod_{i=1}^I \Big(\frac{\lambda_i}{\tilde \lambda_i}\Big)^{\tilde \lambda_i n/p}\prod_{i=1}^I \Big(\int f_i\, d\gamma_{\Sigma_i}\Big)^{\tilde \lambda_i}.$$
\end{proposition}

\hfill \\

\begin{remark}
Together with \cref{Multiple_BThm}, this proves \cref{Main_Thm_IneqPart} by restricting the statement to characteristic functions. The ``Consequently'' part is by optimization over $\tilde \lambda \in \Delta_I^\circ.$ 
Again, this proposition is a consequence of \cref{Sharp_Constant_BM_Ineq}, as if $$\Gamma \le \prod_{i=1}^I \Big(\frac{\tilde \lambda_i}{\lambda_i}\Big)^{2(\frac1p-1)\tilde \lambda_i}\mathcal O_{\tilde \lambda}(\Sigma_1,\cdots,\Sigma_I),$$ then $$\Big[\frac{\det \Gamma_E}{\prod_{i=1}^I (\det((\Sigma_i)_E))^{\tilde \lambda_i}} \Big]^{1/2} \le \prod_{i=1}^I \Big(\frac{\tilde \lambda_i}{\lambda_i}\Big)^{\dim(E)(\frac1p-1)\tilde \lambda_i},$$ and so, the sharp constant for the inequality $$C\gamma_\Gamma\Big(\sum_{i=1}^I \lambda_i K_i\Big) \ge \prod_{i=1}^I \gamma_{\Sigma_i}(K_i)^{\tilde \lambda_i}$$ is at most $\prod_{i=1}^I (\frac{\tilde \lambda_i}{\lambda_i})^{\tilde \lambda_i n/p}.$ For the purpose of studying equality cases, we employ \cref{Gaussian_Case_MultipleBThm}. 
\end{remark}

\hfill \\

\begin{proof}
The Gaussian FRBL associated to the inequality is $$C_{\boldsymbol{\lambda},\boldsymbol{\tilde \lambda}}=\sup \Big\{ \frac{\det(I_n+\Gamma B)}{\prod_{i=1}^I \det(I_n+\Sigma_i A_i)^{\tilde \lambda_i}}: \mathrm{diag}(\tilde \lambda_1 A_1,\cdots,\tilde \lambda_I A_I) \ge (\lambda_i \lambda_j e^{S_i}Be^{S_j})_{1 \le i,j \le I} \Big\}^{1/2}.$$ If we write $D=\sum_{i=1}^I \tilde \lambda_i \log\frac{\tilde \lambda_i}{\lambda_i} \ge 0$, the goal is to prove that $$C_{\boldsymbol{\lambda},\boldsymbol{\tilde \lambda}} \le e^{nD/p}.$$ 
\hfill \\
Define $c=e^{D/p}$ and we want to show that $$C_{\boldsymbol{\lambda},\boldsymbol{\tilde \lambda}} \le c^n.$$  Set $$R_i=S_i+\log \Big(c\frac{\lambda_i}{\tilde \lambda_i}\Big)I_n,$$ so that $$e^{R_i}=c \Big(\frac{\tilde \lambda_i}{\lambda_i}\Big)^{-1}e^{S_i}\quad \forall 1 \le i \le I.$$
\hfill \\
Therefore, using the homogeneity of the Karcher mean, $$\begin{array}{lll}
\ds \Gamma &\le& \ds  \prod_{i=1}^I \Big(\frac{\tilde \lambda_i}{\lambda_i}\Big)^{2(\frac1p-1)\tilde \lambda_i}\mathcal O_{\tilde \lambda }(e^{S_1}\Sigma_1 e^{S_1},\cdots,e^{S_I}\Sigma_I e^{S_I}) \\ 
\nm
&=& \ds c^2\prod_{i=1}^I \Big(\frac{\tilde \lambda_i}{\lambda_i}\Big)^{-2\tilde \lambda_i}\mathcal O_{\tilde \lambda}(e^{S_1}\Sigma_1e^{S_1},\cdots,e^{S_I}\Sigma_Ie^{S_I})=\mathcal O_{\tilde \lambda}(e^{R_1}\Sigma_1 e^{R_1},\cdots,e^{R_I}\Sigma_Ie^{R_I}).
\end{array}$$ 

\hfill \\
Hence, by \cref{Gaussian_Case_MultipleBThm}, as $$\mathrm{diag}(\tilde \lambda_1 A_1,\cdots,\tilde \lambda_I A_I) \ge \Big(\tilde \lambda_i \tilde \lambda_j e^{R_i}\frac B{c^2}e^{R_j} \Big)_{1 \le i,j \le I},$$ we get that $$\prod_{i=1}^I \det(I_n+\Sigma_i A_i)^{\tilde \lambda_i} \ge \det\Big(I_n+\Gamma \frac B{c^2}\Big).$$ Consequently, if $b_1,\cdots,b_n$ denote the eigenvalues of $\Gamma^{1/2}B\Gamma^{1/2}$, $$\begin{array}{lll}
\displaystyle \Bigg( \frac{\det(I_n+\Gamma B)}{\prod_{i=1}^I \det(I_n+\Sigma_i A_i)^{\tilde \lambda_i}}\Bigg)^{1/2} &\le& \displaystyle  \Bigg(\frac{\det(I_n+\Gamma B)}{\det(I_n+c^{-2}\Gamma B)}\Bigg)^{1/2} \\ \\
&\le& \displaystyle c^n \prod_{k=1}^n \Big( \frac{1+b_k}{c^2+b_k}\Big)^{1/2} \le c^n,
\end{array}$$ as $c \ge 1.$ 
\end{proof}

\hfill \\
Thanks to this, we may prove the following Borell-Brascamp-Lieb (BBL) inequality.

\begin{theorem} \label{BBL_Gaussian}
Set $p \in [0;1], \lambda \in \Delta_I^\circ$ and consider $S_1,\cdots,S_I \in \mathrm{Sym}(n), \Gamma, \Sigma_1,\cdots,\Sigma_I \in \mathrm{Sym}_{>0}(n)$ such that $$\Gamma \le \prod_{i=1}^I \Big(\frac{\tilde \lambda_i}{\lambda_i}\Big)^{2(\frac1p-1)\tilde \lambda_i}\mathcal O_{\tilde \lambda}(e^{S_1}\Sigma_1e^{S_1},\cdots,e^{S_I}\Sigma_Ie^{S_I}) \quad \forall \tilde \lambda \in \Delta_I^\circ \text{ if } p \in (0;1],$$ and $$\Gamma \le \mathcal O_\lambda(e^{S_1}\Sigma_1e^{S_1},\cdots,e^{S_I}\Sigma_Ie^{S_I}) \text{ if } p=0.$$  
If $f_i$ are log-concave, centered with respect to $\gamma_{\Sigma_i}$ for any $1 \le i \le I$ and $h$ log-concave such that $$h\Big(\sum_{i=1}^I \lambda_i e^{S_i}x_i \Big) \ge \mathcal M_{\tilde p}^\lambda (f_1(x_1),\cdots,f_I(x_I)) \quad \forall x_1,\cdots,x_I \in \R^n,$$ implies \begin{equation} \label{BBL_Ineq} \int_{\R^n} h \, d\gamma_\Gamma \ge \mathcal M_{\frac{p\tilde p}{p+n\tilde p}}^\lambda\Big(\int_{\R^n}f_1\, d\gamma_{\Sigma_1}, \cdots, \int_{\R^n}f_I\, d\gamma_{\Sigma_I} \Big),\end{equation} for any $\tilde p \ge 0.$ 
If $\tilde p=0$ or $p=0$, $\frac{p\tilde p}{p+n\tilde p}$ is understood to be equal to $0$, while if $\tilde p=+\infty$, $\frac{p\tilde p}{p+n\tilde p}$ is understood to be equal to $\frac pn.$
\end{theorem}

\begin{proof}
If $p=0$ or $\tilde p=0$, it suffices to apply the Forward-Reverse Brascamp-Lieb \cref{FRBL_MultipleBthm} of \cref{Multiple_BThm}, using $\mathcal M_s \ge \mathcal M_0$ for $s \ge 0.$

\hfill \\
With $p,\tilde p>0$ finite, assume that $$h\Big(\sum_{i=1}^I \lambda_i e^{S_i}x_i \Big) \ge \mathcal M_{\tilde p}^\lambda (f_1(x_1),\cdots,f_I(x_I)) \ge \prod_{i=1}^I \Big(\frac{\lambda_i}{\tilde \lambda_i}\Big)^{\tilde \lambda_i/\tilde p}\prod_{i=1}^I f_i(x_i)^{\tilde \lambda_i}.$$ Then, applying the previous result, we get $$\int h \, d\gamma_\Gamma \ge \prod_{i=1}^I \Big(\frac{\lambda_i}{\tilde \lambda_i}\Big)^{\tilde \lambda_i(\frac np+\frac1{\tilde p})} \prod_{i=1}^I \Big(\int_{\R^n} f_i \, d\gamma_{\Sigma_i} \Big)^{\tilde \lambda_i},$$ so that setting $$\frac1r=\frac np+\frac1{\tilde p} \iff r=\frac{p\tilde p}{p+n \tilde p},$$ we get $$\Big(\int_{\R^n} h \, d\gamma_\Gamma \Big)^r \ge \prod_{i=1}^I \bigg[ \frac{\lambda_i}{\tilde \lambda_i} \Big(\int_{\R^n} f_i \, d\gamma_{\Sigma_i} \Big)^r \bigg]^{\tilde \lambda_i}.$$ Choosing the optimal $\tilde \lambda$, we set $$\tilde \lambda_i:=\frac{\lambda_i \Big(\int_{\R^n}f_i \, d\gamma_{\Sigma_i} \Big)^r}{\sum_{j=1}^I \lambda_j \Big(\int_{\R^n} f_j \, d\gamma_{\Sigma_j} \Big)^r} \text{ for any } 1 \le i \le I,$$ and this implies that $$\Big(\int_{\R^n}h \, d\gamma_\Gamma \Big)^r \ge \sum_{i=1}^I \lambda_i \Big(\int_{\R^n}f_i \, d\gamma_{\Sigma_i} \Big)^r.$$

\hfill \\
Finally, with $\tilde p=+\infty$, we may use $$\max_{1 \le i \le I}f_i \ge \prod_{i=1}^I f_i^{\tilde \lambda_i},$$ where $$\tilde \lambda_i:=\frac{\lambda_i \Big(\int_{\R^n}f_i \, d\gamma_{\Sigma_i} \Big)^r}{\sum_{j=1}^I \lambda_j \Big(\int_{\R^n} f_j \, d\gamma_{\Sigma_j} \Big)^r} \text{ for any } 1 \le i \le I,$$ and apply \cref{Computation_towardsBBL}. 
\end{proof}

\hfill \\
In the already interesting special case $\Sigma_i=\Gamma=I_n$ for any $1 \le i \le I$, the Karcher condition may be turned into a more familiar condition thanks to \cref{Relaxed_Conditions_LpBM}. 

\begin{corollary} \label{BBL_Identity_Cov}
    Set $p \in [0;1], \lambda \in \Delta_I^\circ$. Consider $S_1,\cdots,S_I \in \mathrm{Sym}(n)$ such that, if $p \in (0;1)$, $$\Big(\sum_{i=1}^I \lambda_ie^{qS_i}\Big)^{1/q} \ge I_n,$$ where $q=\frac p{p-1} \in (-\infty;0)$, while if $p=0$, $$\sum_{i=1}^I \lambda_i S_i \ge 0,$$ and if $p=1$, $$S_i \ge 0 \quad \forall 1 \le i \le I.$$
    Then, if $f_i$ are log-concave, centered with respect to the standard Gaussian measure $\gamma$ and $h$ log-concave, $$h\Big(\sum_{i=1}^I \lambda_i e^{S_i}x_i\Big) \ge \mathcal M_{\tilde p}^\lambda(f_1(x_1),\cdots,f_I(x_I)) \quad \forall x_1,\cdots,x_I \in \R^n,$$ implies $$\int_{\R^n}h\, d\gamma_n \ge \mathcal M_{\frac{p\tilde p}{p+n\tilde p}}^\lambda\Big(\int_{\R^n}f_1\, d\gamma_n,\cdots,\int_{\R^n}f_I\, d\gamma_n\Big),$$ for any $\tilde p \ge 0.$  
\end{corollary}

\hfill \\
Let us now investigate the equality cases of Theorem \ref{BBL_Gaussian}, that was obtained as optimization of a family of FRBL. In particular, if we have equality in BBL, this means that we have equality in a certain ``optimal'' FRBL, which allows to run a similar argument, to the equality cases of Theorem \ref{Multiple_BThm}. Let us restrict to the case $p,\tilde p>0.$ 

\hfill \\
More precisely, assume $$\Gamma \le \prod_{i=1}^I \Big(\frac{\tilde \lambda_i}{\lambda_i}\Big)^{2(\frac1p-1)\tilde \lambda_i}\mathcal O_{\tilde \lambda}(e^{S_1}\Sigma_1e^{S_1},\cdots,e^{S_I}\Sigma_Ie^{S_I}) \quad \forall \tilde \lambda \in \Delta_I^\circ.$$ Further, suppose that $f_i$ are log-concave, centered with respect to $\gamma_{\Sigma_i}$ for any $1 \le i \le I$ and $h$ log-concave, with $$h\Big(\sum_{i=1}^I \lambda_i e^{S_i}x_i\Big) \ge \mathcal M_{\tilde p}^\lambda(f_1(x_1),\cdots,f_I(x_I)),$$ for any $x_1,\cdots,x_I \in \R^n$ and satisfying equality in the BBL: $$\Big(\int_{\R^n}h\, d\gamma_\Gamma \Big)^r=\sum_{i=1}^I \lambda_i \Big(\int_{\R^n}f_i\, d\gamma_{\Sigma_i} \Big)^r,$$ where $$r=\frac{p\tilde p}{p+n\tilde p}.$$
\hfill \\
This means that $((f_i)_{1 \le i \le I},\tilde h)$ satisfies exactly equality for the FRBL of \cref{Computation_towardsBBL}, \textit{i.e.} setting $$\tilde \lambda_i=\frac{\lambda_i \Big(\int_{\R^n}f_i\, d\gamma_{\Sigma_i}\Big)^r}{\sum_{j=1}^I \lambda_j \Big(\int_{\R^n}f_j\, d\gamma_{\Sigma_j} \Big)^r} \text{ for any } 1 \le i \le I,$$ with $$\prod_{i=1}^I \Big(\frac{\lambda_i}{\tilde \lambda_i}\Big)^{-\tilde \lambda_i/\tilde p}h\Big(\sum_{i=1}^I \lambda_i e^{S_i}x_i\Big):=\tilde h\Big(\sum_{i=1}^I \lambda_i e^{S_i}x_i\Big) \ge \prod_{i=1}^I f_i(x_i)^{\tilde \lambda_i},$$ one has $$\prod_{i=1}^I \Big(\frac{\tilde \lambda_i}{ \lambda_i}\Big)^{\tilde \lambda_i n/p}\int_{\R^n} \tilde h\, d\gamma_\Gamma=\prod_{i=1}^I \Big(\int_{\R^n}f_i\, d\gamma_{\Sigma_i} \Big)^{\tilde \lambda_i}.$$  

\hfill \\
Note that if $\tilde \lambda_i=\lambda_i$, we already know that the Gaussian FRBL constant is exactly $1$, under the assumption of $$\Gamma \le \mathcal O_\lambda(e^{S_1}\Sigma_1e^{S_1},\cdots,e^{S_I}\Sigma_Ie^{S_I}),$$  by \cref{Multiple_BThm}.
The assumption in our context is clearly stronger than the one of \cref{Multiple_BThm} (take $\tilde \lambda=\lambda$) and we get the desired equality cases by \cref{Eq_Cases_MultipleBThm}. 

\hfill \\
If $\tilde \lambda_i \not =\lambda_i$ for some $1 \le i \le I$, this implies that the Gaussian FRBL constant is equal to $$\prod_{i=1}^I \Big(\frac{\tilde \lambda_i}{\lambda_i}\Big)^{\tilde \lambda_i n/p}=c^n,$$ where $c=e^{D/p}>1$ with $$D=\sum_{i=1}^I \tilde \lambda_i \log \frac{\tilde \lambda_i}{\lambda_i}>0.$$

\hfill \\
Indeed, from the proof of \cref{Computation_towardsBBL}, it is enough to notice that if the Gaussian FRBL constant is equal to $c^n$, then it cannot be attained by any Gaussian matrices. Indeed, for $b_1,\cdots,b_n \ge 0$ eigenvalues of $\Gamma^{1/2}B\Gamma^{1/2}$, $$\prod_{k=1}^n \Big(\frac{1+b_k}{c^2+b_k}\Big)^{1/2}<1.$$

\hfill \\
Let us show that this implies that there cannot be any extremizers in this case.
Indeed, by a similar doubling argument to the one already explained in the second section, if there were extremizers, the constant would be attained at certain positive semi-definite Gaussian matrices, which is not the case by what we just proved. 

\begin{lemma}
    Assume that there exists a finite constant $C<+\infty$ such that for any $f_i$ log-concave, centered with respect to $\gamma_{\Sigma_i}$ for any $1 \le i \le I$ and $h$ log-concave, $$h\Big(\sum_{i=1}^I \lambda_i e^{S_i}x_i\Big) \ge \prod_{i=1}^I f_i(x_i)^{\tilde \lambda_i} \Longrightarrow \int_{\mathbb{R}^n} h \, d\gamma_\Gamma \ge C^{-1}\prod_{i=1}^I \Big(\int_{\mathbb{R}^n}f_i\, d\gamma_{\Sigma_i}\Big)^{\tilde \lambda_i}.$$ Then, if equality is attained by a tuple with positive finite integrals, $C$ is attained by a tuple of finite positive semi-definite matrices. 
\end{lemma}

\begin{proof}
    After normalization, we may write the extremizing pointwise inequality as $$H\Big(\sum_{i=1}^I \lambda_i e^{S_i}x_i\Big) \ge C\prod_{i=1}^I F_i(x_i)^{\tilde \lambda_i}, \quad \int_{\mathbb{R}^n} H\, d\gamma_\Gamma=\int_{\mathbb{R}^n} F_i\, d\gamma_{\Sigma_i}=1.$$ Then, the normalized iterates satisfy $$\mathrm{Conv}^k_{\gamma_\Gamma}(H)\Big(\sum_{i=1}^I \lambda_i e^{S_i}x_i\Big) \ge C\prod_{i=1}^I \mathrm{Conv}_{\gamma_{\Sigma_i}}^k(F_i)(x_i)^{\tilde \lambda_i},$$ for every $k \ge 1$, by using the doubling step and then applying the Forward-Reverse Brascamp-Lieb inequality. One may then again show that $H$ must be centered, so that we can apply the local central limit theorem with normalized Gaussian densities $G_i$ and $G_0$ satisfying the pointwise inequality. Let $A_i=\mathrm{Cov}(F_i\gamma_{\Sigma_i})^{-1}-\Sigma_i^{-1}\ge 0$ and $A_0=\mathrm{Cov}(H\gamma_{\Gamma})^{-1}-\Gamma^{-1} \ge 0.$ The limiting pointwise inequality gives the block constraint: $$\mathrm{diag}(\tilde \lambda_1 A_1,\cdots,\tilde \lambda_I A_I) \ge (\lambda_i \lambda_j e^{S_i}A_0e^{S_j})_{1 \le i,j \le I},$$ while evaluating the inequality at $x_i=0$ gives $$\Big(\frac{\det(I_n+\Gamma A_0)}{\prod_{i=1}^I \det(I_n+\Sigma_iA_i)^{\tilde \lambda_i}}\Big)^{1/2} \ge C.$$  Moreover, the Gaussian bound is at most $C$, by applying the Forward-Reverse Brascamp-Lieb inequality to Gaussians. Thus, equality for Gaussians holds as desired. 
\end{proof}

\hfill \\
This allows to assemble the inequality and the equality cases in the following unified theorem. 

\begin{theorem} \label{Complete_BBL_eq}
Set $p \in [0;1], \tilde p \in [0;+\infty], \lambda \in \Delta_I^\circ$ and $S_1,\cdots,S_I \in \mathrm{Sym}(n), \Gamma, \Sigma_1,\cdots,\Sigma_I \in \mathrm{Sym}_{>0}(n)$ such that $$\Gamma \le \prod_{i=1}^I \Big(\frac{\tilde \lambda_i}{\lambda_i}\Big)^{2(\frac1p-1)\tilde \lambda_i}\mathcal O_{\tilde \lambda}(e^{S_1}\Sigma_1e^{S_1},\cdots,e^{S_I}\Sigma_Ie^{S_I}) \quad \forall \tilde \lambda \in \Delta_I^\circ \text{ if } p \in (0;1],$$ and $$\Gamma \le \mathcal O_\lambda(e^{S_1}\Sigma_1e^{S_1},\cdots,e^{S_I}\Sigma_Ie^{S_I}) \text{ if } p=0.$$  
Define $$E=\bigcap_{i=1}^I \ker(e^{S_i}\Sigma_i e^{S_i}-\Gamma).$$
Assume that $f_i$ is log-concave and centered with respect to $\gamma_{\Sigma_i}$ for any $1 \le i \le I$ and $h$ log-concave, such that $$0<\int_{\mathbb{R}^n}f_i\, d\gamma_{\Sigma_i}, \; \int_{\mathbb{R}^n}h\, d\gamma_\Gamma<+\infty,$$ and $$h\Big(\sum_{i=1}^I \lambda_i e^{S_i}x_i \Big) \ge \mathcal M_{\tilde p}^\lambda (f_1(x_1),\cdots,f_I(x_I)) \quad \forall x_1,\cdots,x_I \in \R^n.$$ 
Then, $$\int_{\mathbb{R}^n} h \, d\gamma_\Gamma \ge \mathcal M_{\frac{p\tilde p}{p+n\tilde p}}^\lambda \Big(\int_{\mathbb{R}^n}f_1\, d\gamma_{\Sigma_1},\cdots,\int_{\mathbb{R}^n}f_I\, d\gamma_{\Sigma_I}\Big),$$ with the common endpoint conventions introduced in \cref{BBL_Gaussian}. 
The equality cases depend on the pointwise exponent $\tilde p$:
\begin{itemize}
\item[(i)] If $\tilde p=0$, there is equality iff there exist constants $a_1,\cdots,a_I>0$ and a non-zero log-concave function $\phi: E \longrightarrow [0;+\infty)$, centered w.r.t $\gamma_{\Gamma_E}$, so that $$f_i(x)=a_i \phi(P_{E}e^{S_i}x) \quad \forall 1 \le i \le I, h(x)=\prod_{i=1}^I a_i^{\lambda_i} \phi(P_Ex),$$ for a.e. $x \in \R^n.$ 
\item[(ii)] If $0<\tilde p<+\infty$, equality holds exactly when there are a single constant $a>0$ and a non-zero centered $\phi$ on $E$ such that $\phi^{\tilde p}$ is concave on its convex positive support and $$f_i(x)=a \phi(P_Ee^{S_i}x), 
\quad h=a\phi( P_Ex),$$ for a.e. $x \in \R^n.$
\item[(iii)] If $\tilde p=+\infty$, equality holds exactly when $$f_i=a\mathbbm{1}_{e^{-S_i}(C+E^\perp)}, \quad h=a\mathbbm{1}_{C+E^\perp} \text{ almost everywhere},$$ for $a>0$ and a centered convex set $C \subset E$ of positive Gaussian measure. 
\end{itemize}
\end{theorem}

\begin{proof}
    First, (i) is an immediate consequence of \cref{Eq_Cases_MultipleBThm}. If $p>0$ and $\tilde p>0$, the preceding saturation argument forces $$\int f_1\, d\gamma_{\Sigma_1}=\cdots=\int f_I\, d\gamma_{\Sigma_I}.$$ The BBL inequality then becomes \cref{Eq_Cases_MultipleBThm}, so that there is a common profile $\phi$ as in (i). The equality of integrals of $f_i$ forces all the scalar factors to agree. 
    \hfill \\
    If $p=0$ and $\tilde p>0$, equality in \cref{BBL_Ineq} again forces equality in 
    \cref{Eq_Cases_MultipleBThm}, so that we have proportional common profiles. Evaluating the stronger pointwise assumption at a point where the profile is positive yields that $$\prod_{i=1}^I a_i^{\lambda_i} \ge \mathcal M_{\tilde p}^\lambda(a_1,\cdots,a_I) \Longrightarrow a_1=\cdots=a_I.$$ Moreover, canceling the scalar contributions, the pointwise assumption becomes $$\phi\Big(\sum_{i=1}^I \lambda_i x_i\Big) \ge \mathcal M_{\tilde p}^\lambda(\phi(x_1),\cdots,\phi(x_I)),$$ on its positive support, which implies concavity of $\phi^{\tilde p}$ and concludes the necessity of (ii). Indeed, continuity on the relative interior of the support implies that the Jensen concavity only for fixed weights is upgraded to ordinary concavity. 
    \hfill \\
    If $\tilde p=+\infty$, we find again that $$f_i=a\phi(P_Ee^{S_i}\cdot), \quad h=a\phi \circ P_E.$$ Then, the pointwise assumption now says that $$\phi\Big(\sum_{i=1}^I \lambda_i x_i\Big) \ge \max_{1 \le i \le I}\phi(x_i) \text{ whenever } \phi(x_i)>0.$$ Let $U=\mathrm{int}_E \{\phi>0\}.$ We claim that $\phi$ is constant on $U.$ Indeed, for $x,y \in U$ and a fixed $1 \le j \le I$, we have that $$\phi(\lambda_j x+(1-\lambda_j)y) \ge \max\{\phi(x),\phi(y)\},$$ so that defining $x_0=x, x_{k+1}=\lambda_j x_k+(1-\lambda_j)y$, all $x_k$ belong to $U$ by convexity and $x_k=\lambda_j^kx+(1-\lambda_j^k)y \rightarrow_{k \rightarrow +\infty} y.$ As $\phi$ is continuous on $U$ and $\phi(x_k) \ge \phi(x)$ for any $k \ge 1$, we get that $\phi(y) \ge \phi(x)$ and we may exchange the roles of $x,y$ to obtain the desired claim. Consequently, up to a null set, $$\phi=c \mathbbm{1}_C \text{ for some convex set } C \subset E \text{ with }\mathrm{bar}_{\gamma_{\Gamma_E}}(C)=0.$$ 
    \hfill \\
    The sufficiency is quite clear in all cases. 
\end{proof}

\section{From centered to uncentered}

In this section, we provide a new slicing approach to be able to extend a Gaussian inequality applied to centered, convex sets to sets with possibly arbitrary barycenters. 

\subsection{From the point of view of inequalities}

Consider $N \ge 2$ and $$H_N=\left \{ z=(z_1,\cdots,z_N) \in (\R^n)^N : \sum_{j=1}^N z_j=0 \right\},$$ a Euclidean space of dimension $n(N-1).$ For $K \subset \R^n$ and $a \in \R^n$, define $$C_N(K,a)=H_N \cap (K-a)^N.$$ 
\hfill \\
For $\Sigma \in \mathrm{Sym}_{>0}(n)$, we write $\gamma_{\Sigma,H_N}$ as the Gaussian probability measure on $H_N$ with covariance $$\Sigma_{H_N}=(I_N \otimes \Sigma)\big|_{H_N}.$$ As $H_N \backsimeq \R^{N-1} \otimes \R^n$, we may identify this covariance with $I_{N-1} \otimes \Sigma.$ 

\hfill \\
Observe that $C_N(K,a)$ satisfies the following properties: 
\begin{itemize}
    \item [(i)] If $a \in \mathrm{int}(K)$ , the set $C_N(K,a)$ has positive $\gamma_{H_N}$-measure. Indeed, since $a \in \mathrm{int}(K)$, a small Euclidean ball around the origin is contained in $K-a$ and its $N$-times product intersects $H_N$ in a neighborhood of the origin.  
    \item [(ii)] If $K$ is convex, $C_N(K,a)$ is also convex.
    \item[(iii)] The barycenter of $C_N(K,a)$ w.r.t $\gamma_{H_N}$ is zero whenever the slice has positive measure, or in other words, $$\int_{C_N(K,a)}z\, d\gamma_{\Sigma,H_N}(z)=0.$$ Indeed, $H_N, \gamma_{H_N}$ and $C_N(K,a)$ are invariant under permutation of the coordinates, which means that the Gaussian barycenter of $C_N(K,a)$ is also invariant under such permutations. In particular, it needs to belong to $$H_N \cap \{ (v,\cdots,v) \in (\R^n)^N\}=\{0\}.$$
    \hfill \\
    The same argument applies to $\gamma_{\Sigma,H_N}$ and to the convex hull of $C_N(K,a).$ The centering of this convex hull follows again from the permutation invariance. 
    \item[(iv)] Since $H_N$ is convex, $$\sum_{i=1}^I \lambda_i C_N(K_i,a_i) \subset C_N \Big(\sum_{i=1}^I \lambda_i K_i, \sum_{i=1}^I \lambda_i a_i \Big). $$
\end{itemize}
\hfill \\
Let $K \subset \R^n$ be measurable with $\gamma_\Sigma(K)>0$ and define $$\gamma_{\Sigma,K}=\frac{\mathbbm{1}_K \gamma_\Sigma}{\gamma_\Sigma(K)}.$$
For $a \in \mathbb{R}^n$, set $f_{N,a}$ the density of $\frac1{\sqrt N}\sum_{j=1}^N (X_j-a)$, where $X_1,\cdots,X_N$ i.i.d, with law $\gamma_{\Sigma,K}$, or in other words, $$f_{N,a}(y)=N^{n/2}(\gamma_{\Sigma,K})^{*,N}(Na+\sqrt{N}y) \quad \forall y \in \mathbb{R}^n.$$

\begin{lemma} \label{High_Dim_Slicing}
    Let $K$ be measurable with $\gamma_\Sigma(K)>0.$ Then, for every $a \in \R^n$ and $N \ge 2$, $$\gamma_{\Sigma,H_N}(C_N(K,a))=(2\pi)^{n/2}(\det \Sigma)^{1/2}e^{\frac N2 \langle \Sigma^{-1}a,a\rangle}\gamma_\Sigma(K)^Nf_{N,a}(0).$$  
\end{lemma}

\begin{proof}
    For every $a \in \mathbb{R}^n$ and $w \in H_N$, $$\sum_{j=1}^N \langle \Sigma^{-1}(a+w_j),a+w_j\rangle=N\langle \Sigma^{-1}a,a\rangle+\sum_{j=1}^N \langle \Sigma^{-1}w_j,w_j\rangle,$$ as $$2\Big\langle \Sigma^{-1}a,\sum_{j=1}^N w_j\Big\rangle=0.$$ Hence, \begin{equation} \label{Gaussian_Factorization} \prod_{j=1}^N \gamma_\Sigma(a+w_j)=\gamma_{\Sigma}(\sqrt{N}a)\gamma_{\Sigma,H_N}(w),\end{equation} as $$\gamma_{\Sigma,H_N}(w)=\frac{\exp\Big[-\frac12 \sum_{j=1}^N \langle \Sigma^{-1}w_j,w_j\rangle \Big]}{(2\pi)^{n(N-1)/2}(\det \Sigma)^{(N-1)/2}},$$ and by abuse of notation, we write $$\gamma_\Sigma(a)=\frac1{(2\pi)^{n/2}(\det \Sigma)^{1/2}}e^{-\frac12 \langle \Sigma^{-1}a,a\rangle}.$$In addition, we claim that $$\int_{H_N}\prod_{j=1}^N \gamma_{\Sigma,K}(a+w_j)\, d_Hw=f_{N,a}(0).$$ For fixed $a$, consider $x_1,\cdots,x_{N-1} \in \mathbb{R}^n$ and set $$w_j=x_j-a \text{ for } j<N, \quad w_N=(N-1)a-\sum_{j=1}^{N-1}x_j.$$ Then, $\sum_j w_j=0$ and every $w \in H_N$ may be obtained uniquely in this way. The differential of the parametrization is the linear map $$L(h_1,\cdots,h_{N-1})=\Big(h_1,\cdots,h_{N-1},-\sum_{j=1}^{N-1}h_j\Big),$$ so that its Gram matrix is $$L^*L=\Big(I_{N-1}+\mathbbm{1} \otimes \mathbbm{1}\Big) \otimes I_n,$$ where $\mathbbm{1}=(1,\cdots,1) \in \mathbb{R}^{N-1}.$ Therefore, $$\det(L^*L)=N^n,$$ and so, $$d_Hw=N^{n/2}dx_1\cdots dx_{N-1}.$$ Hence, $$\begin{array}{lll}
    \ds \int_{H_N}\prod_{j=1}^N \gamma_{\Sigma,K}(a+w_j)\, d_Hw &=& \ds N^{n/2}\int_{(\mathbb{R}^n)^{N-1}}\prod_{j=1}^{N-1}\gamma_{\Sigma,K}(x_j) \cdot \gamma_{\Sigma,K}\Big(Na-\sum_{j=1}^{N-1}x_j\Big)\, dx_1\cdots dx_{N-1}
    \\
    \nm 
    &=& \ds N^{n/2}(\gamma_{\Sigma,K})^{*N}(Na)=f_{N,a}(0).
    \end{array}$$
    Multiplying \cref{Gaussian_Factorization} by $\prod_{j=1}^N \mathbbm{1}_K(a+w_j)$ and integrating over $H_N$, $$\begin{array}{lll}
    \ds \gamma_\Sigma(\sqrt N a)\gamma_{\Sigma,H_N}(C_N(K,a))  &=& \ds \gamma_\Sigma(\sqrt N a)\int_{H_N}\Big(\prod_{j=1}^N \mathbbm{1}_K(a+w_j)\Big)\gamma_{\Sigma,H_N}(w)\, d_Hw \\
    \nm 
    &=& \ds \int_{H_N}\prod_{j=1}^N \Big[ \mathbbm{1}_K(a+w_j) \gamma_\Sigma(a+w_j) \Big]\, d_Hw \\
    \nm 
    &=& \ds \gamma_\Sigma(K)^N \int_{H_N}\prod_{j=1}^N \gamma_{\Sigma,K}(a+w_j) \, d_Hw=\gamma_\Sigma(K)^N f_{N,a}(0).
    \end{array}$$ Dividing by $\gamma_{\Sigma}(\sqrt Na)$ and using its definition, we conclude the proof of the lemma. 
\end{proof}

\begin{remark}
We now study the behavior of $f_{N,a}$ as $N \rightarrow +\infty.$ 
\hfill \\
First, we claim that $$f_{N,a}(0) \le \frac{N^{n/2}}{(2\pi)^{n/2}(\det \Sigma)^{1/2}\gamma_\Sigma(K)}.$$ Indeed, if $X \sim \gamma_{\Sigma,K}$, $$\|\gamma_{\Sigma,K}\|_{L^\infty} \le \frac1{(2\pi)^{n/2}(\det \Sigma)^{1/2}\gamma_\Sigma(K)}.$$ Thus, $$\|f_{N,a}\|_{L^\infty}=N^{n/2}\|\gamma_{\Sigma,K}^{*N}\|_{L^\infty} \le N^{n/2}\|\gamma_{\Sigma,K}\|_{L^\infty} \le \frac{N^{n/2}}{(2\pi)^{n/2}(\det \Sigma)^{1/2}\gamma_\Sigma(K)}.$$ In particular, $$\gamma_{\Sigma,H_N}(C_N(K,a)) \le e^{N\langle \Sigma^{-1}a,a\rangle/2} N^{n/2} \gamma_\Sigma(K)^{N-1},$$ so that \begin{equation} \label{High_Dim_Gaussian_Limsup} \limsup_{N \rightarrow +\infty}\gamma_{\Sigma,H_N}(C_N(K,a))^{1/(N-1)n} \le e^{\langle \Sigma^{-1}a,a\rangle/2n}\gamma_\Sigma(K)^{1/n}.\end{equation}
\hfill \\
For $K$ an arbitrary measurable set with positive Gaussian measure, we can employ the local central limit theorem to obtain that $$f_{N,b_{K,\Sigma}}(0) \rightarrow_{N \rightarrow +\infty} \frac1{(2\pi)^{n/2}\sqrt{\det(\mathrm{Cov}(\gamma_{\Sigma,K}))}}>0.$$
In particular, $$\lim_{N \rightarrow +\infty} f_{N,b_{K,\Sigma}}^{1/N}(0)=1.$$
\end{remark}

\hfill \\
By \cref{High_Dim_Slicing} and the previous remark, \begin{equation} \label{High_Dim_Gaussian_Lim} \lim_{N \rightarrow +\infty} \gamma_{\Sigma,H_N}(C_N(K,b_{K,\Sigma}))^{1/(N-1)n}=e^{\langle \Sigma^{-1}b_{K,\Sigma},b_{K,\Sigma}\rangle/2n}\gamma_\Sigma(K)^{1/n}.\end{equation}

\hfill \\
We may now prove the generalized uncentered Gaussian Brunn-Minkowski. 

\begin{proof}[Proof of \cref{Generalized_Gaussian_BM}]
It suffices to apply the Gaussian Brunn-Minkowski theorem, which is the case $p=1,\tilde p=\infty \text{ and } S_i=0$ of \cref{BBL_Identity_Cov}, on the higher-dimensional space $H_N$ with the sets $$\{\mathrm{conv}(C_N(K_i,b_{K_i}))\}_{1 \le i \le I},$$ which are the convex hulls of $C_N(K_i,b_{K_i})$, and then let $N \rightarrow +\infty$, using \eqref{High_Dim_Gaussian_Limsup} and \eqref{High_Dim_Gaussian_Lim}.
\hfill \\
Indeed, observe that $$\sum_{i=1}^I \lambda _i \mathrm{conv}(C_N(K_i,b_{K_i}))=\mathrm{conv}\Big(\sum_{i=1}^I \lambda_i C_N(K_i,b_{K_i}) \Big)   \subset C_N \Big(\sum_{i=1}^I \lambda_iK_i,\sum_{i=1}^I \lambda_i b_{K_i}\Big),$$ as $\sum_{i=1}^I \lambda_iK_i$ is convex, which means that $$C_N\Big(\sum_{i=1}^I \lambda_i K_i,\sum_{i=1}^I \lambda_i b_{K_i}\Big) \text{ is convex and contains } \sum_{i=1}^I \lambda_i C_N(K_i,b_{K_i}).$$  
\hfill \\
Then, $$\gamma_{H_N}\Big( C_N\Big(\sum_{i=1}^I \lambda_i K_i, \sum_{i=1}^I \lambda_i b_{K_i}\Big) \Big)^{1/(N-1)n} \ge \sum_{i=1}^I \lambda_i \gamma_{H_N}(C_N(K_i,b_{K_i}))^{1/(N-1)n},$$ and we may let $N \rightarrow +\infty.$ 
The conclusion with Gaussian measure $\gamma_\Sigma$ is a consequence of considering instead $\Sigma^{-1/2}K_i$ for any $1 \le i \le I.$ 
\end{proof}

\hfill \\
A similar strategy may also be applied to Gaussian correlation  (or also strong Gaussian correlation). Regarding such inequalities, it is important for $K_1$ and $K_2$ to have the same Gaussian barycenters since the intersection appears and $$C_N(K_1,b) \cap C_N(K_2,b)=C_N(K_1 \cap K_2,b).$$ Hence, if $K_1,K_2 \subset \R^n$ are convex sets with $\mathrm{bar}_\gamma(K_1)=\mathrm{bar}_\gamma(K_2):=b$, $$\begin{array}{lll}
    \displaystyle \gamma_{H_N}(C_N(K_1 \cap K_2,b))^{1/(N-1)} &=& \gamma_{H_N}(C_N(K_1,b) \cap C_N(K_2,b) )^{1/(N-1)} \\ \\
    &\ge& \gamma_{H_N}(C_N(K_1,b))^{1/(N-1)} \gamma_{H_N}(C_N(K_2,b))^{1/(N-1)},
    \end{array}$$ by centered Gaussian correlation. Letting $N \rightarrow +\infty$, we get the following corollary.
    \begin{corollary}
    If $K_1,K_2 \subset \R^n$ are convex sets, with the same Gaussian barycenter $b$, then
    $$\gamma_n(K_1 \cap K_2) \ge e^{|b|^2/2}\gamma_n(K_1)\gamma_n(K_2).$$
    \end{corollary}
\hfill \\
More generally, the following extension of Gaussian correlation holds. 

\begin{theorem} \label{Uncentered_GCI}
    Let $K_1,\cdots,K_I \subset \R^n$ be convex sets, with $\mathrm{bar}_{\gamma_{\Sigma_i}}(K_i)=b$ for any $1 \le i \le I.$ Set $\Sigma_0^{-1} \ge \Sigma_i^{-1}>0$ for any $1 \le i \le I.$ Then, \begin{equation} \label{Uncentered_GCI_ineq} \gamma_{\Sigma_0}\Big(\bigcap_{i=1}^I K_i\Big) \ge \exp \bigg[\frac12 \Big\langle b, \Big(\sum_{i=1}^I \Sigma_i^{-1}-\Sigma_0^{-1}\Big) b\Big\rangle \bigg]\prod_{i=1}^I \gamma_{\Sigma_i}(K_i).\end{equation} 
\end{theorem}

\begin{proof}
    Again, it suffices to apply \cite[Theorem 1.5]{Shohei_Hiroshi2} to the convex sets $C_N(K_i,b)$, that are centered with respect to $\gamma_{\Sigma_i,H_N}$ and satisfy $$\bigcap_{i=1}^I C_N(K_i,b)=C_N\Big(\bigcap_{i=1}^I K_i,b\Big),$$ and use both \cref{High_Dim_Gaussian_Limsup} and \cref{High_Dim_Gaussian_Lim} to get the desired inequality. 
\end{proof}

\hfill \\

\subsection{From the point of view of equality cases}

\hfill \\
In this subsection, the main goal is to derive the equality cases of Theorem \ref{Generalized_Gaussian_BM} and \cref{Uncentered_GCI}. 

\hfill \\
Since we derived Theorem \ref{Generalized_Gaussian_BM} as a limiting case of the corresponding Gaussian Brunn-Minkowski inequality by adding dimensions, it is not \textit{a priori} clear how to deal with the equality cases. In particular, equality in \cref{Generalized_Gaussian_BM} does not necessarily imply equality in any finite-dimensional inequality with $\{C_N(K_i,b_{k_i})\}_{1 \le i \le I}.$

\hfill \\
Nevertheless, we obtain the following equality cases for the uncentered Gaussian Brunn-Minkowski inequality. 

\begin{theorem}
    Under the hypotheses of \cref{Generalized_Gaussian_BM}, equality in the generalized Gaussian Brunn-Minkowski inequality holds if and only if
    $$\gamma_n(K_i \Delta K)=0 \quad \forall 1 \le i \le I, K=\sum_{i=1}^I \lambda_i K_i.$$
\end{theorem}

\begin{proof}
Trivially, $C_N(K,a)=C_N(K-u,a-u)$ for any $u \in \mathbb{R}^n.$ By the centered Gaussian Brunn-Minkowski, $$\gamma_{H_N}\Big(C_N\Big(\sum_{i=1}^I \lambda_i K_i-u, \sum_{i=1}^I \lambda_i b_{K_i}-u\Big) \Big)^{1/(N-1)n} \ge \sum_{i=1}^I \lambda_i \gamma_{H_N}(C_N(K_i,b_{K_i}))^{1/(N-1)n},$$ for any $u \in \R^n.$ Taking the limit as $N \rightarrow +\infty$ in the right-hand side and the limsup in the left-hand side, $$e^{|\sum_{i=1}^I \lambda_i b_{K_i}-u|^2/2n}\gamma_n\Big(\sum_{i=1}^I \lambda_i K_i-u\Big)^{1/n} \ge \sum_{i=1}^I \lambda_i e^{|b_{K_i}|^2/2n}\gamma_n(K_i)^{1/n},$$ for any $u \in \R^n.$  

\hfill \\
Now, suppose that $K_1,\cdots,K_I$ are measurable sets such that $$e^{|\sum_{i=1}^I \lambda_i b_{K_i}|^2/2n}\gamma_n\Big(\sum_{i=1}^I \lambda_i K_i\Big)^{1/n}=\sum_{i=1}^I \lambda_i e^{|b_{K_i}|^2/2n}\gamma_n(K_i)^{1/n}.$$ 
We may define $$F: \R^n \rightarrow [0;+\infty), \quad F(u)=e^{|\sum_{i=1}^I \lambda_i b_{K_i}-u|^2/2n}\gamma_n\Big(\sum_{i=1}^I \lambda_i K_i-u\Big)^{1/n}.$$ Then, $$0=\nabla \log F(0)=\frac{-\sum_{i=1}^I \lambda_i b_{K_i}+b_{\sum_{i=1}^I \lambda_i K_i}}n,$$ since $$\log \gamma_n(K-u)=\log \int_K (2\pi)^{-n/2}e^{-\frac{|y-u|^2}2}\, dy=-\frac{|u|^2}2+\log \int_K (2\pi)^{-n/2}e^{\langle y,u\rangle}e^{-\frac{|y|^2}2}\, dy,$$ so that differentiating at $u=0$ yields $$\frac{\int_K (2\pi)^{-n/2} ye^{-\frac{|y|^2}2}\, dy}{\int_K (2\pi)^{-n/2} e^{-\frac{|y|^2}2}\, dy}=b_K.$$ Hence, we conclude that $$b_{\sum_{i=1}^I \lambda_i K_i}=\sum_{i=1}^I \lambda_i b_{K_i}.$$ Now, set $r_K(v)=e^{|b_{K+v}|^2/2n}\gamma_n(K+v)^{1/n}$ for $K \subset \R^n, v \in \R^n$ and $$G(v_1,\cdots,v_I)=e^{|\sum_{i=1}^I \lambda_i b_{K_i+v_i}|^2/2n}\gamma_n\Big(\sum_{i=1}^I \lambda_i K_i+\sum_{i=1}^I \lambda_i v_i\Big)^{1/n}-\sum_{i=1}^I \lambda_i r_{K_i}(v_i),$$ for any $v_1,\cdots,v_I \in \R^n.$ 
\hfill \\
We have already seen that $$\nabla \big[ \log \gamma_n(K+v) \big]=-b_{K+v},$$ and since $$b_{K+v}=\frac{\int_{K+v}ye^{-|y|^2/2}\, dy}{\int_{K+v}e^{-|y|^2/2}\, dy} =v+\frac{\int_K ye^{-|y+v|^2/2}\, dy}{\int_K e^{-|y+v|^2/2}\, dy},$$ one has $$D\big(b_{K+v}\big)=I_n-\mathrm{Cov}(\gamma_{K+v}),$$ where $\mathrm{Cov}(\gamma_{K+v})$ is the covariance matrix of the Gaussian measure induced on $K+v$, \textit{i.e.} $$\mathrm{Cov}(\gamma_{K+v})=\frac1{\gamma_n(K+v)}\int_{K+v}x \otimes x\, d\gamma_n(x)-b_{K+v} \otimes b_{K+v}.$$ In particular, $$\nabla \log r_K(0)=\frac1n(I_n-\mathrm{Cov}(\gamma_K))b_K-\frac1n b_K=-\frac1n \mathrm{Cov}(\gamma_K)b_K,$$ and so $$\nabla r_K(0)=-\frac{r_K(0)}n \mathrm{Cov}(\gamma_K)b_K .$$ Note that the quantities considered here are smooth functions of the translation parameter and we may differentiate under the integral sign as the integrands are dominated by Gaussians. Recall that $$G(v_1,\cdots,v_I) \ge 0=G(0,\cdots,0),$$ so that for any $1 \le i \le I$, $$\begin{array}{lll}
 0 &=& \displaystyle \nabla_{v_i}G(0)=-\frac{\lambda_i r_{\sum_{j=1}^I \lambda_j K_j }(0)}{n} \mathrm{Cov}(\gamma_{K_i})b_{\sum_{j=1}^I \lambda_j K_j}+\frac{\lambda_i r_{K_i}(0)}{n}\mathrm{Cov}(\gamma_{K_i})b_{K_i} \\ \\
&=& \displaystyle \frac{\lambda_i}{n}\mathrm{Cov}(\gamma_{K_i}) \Big( r_{K_i}(0)b_{K_i}-r_{\sum_{j=1}^I \lambda_j K_j}(0) b_{\sum_{j=1}^I \lambda_j K_j}\Big),
\end{array}$$ for any $1 \le i \le I.$ 

\hfill \\
Here, we used that $$\begin{array}{lll}
\displaystyle \nabla_{v_i} \Big[ \frac12 \Big| \sum_{j=1}^I \lambda_j b_{K_j+v_j}\Big|^2+\log \gamma_n \Big(\sum_{j=1}^I \lambda_j K_j+\sum_{j=1}^I \lambda_j v_j\Big) \Big]\Big|_{v_1=\cdots=v_I=0}=\lambda_i D(b_{K_i})\Big[ \sum_{j=1}^I \lambda_j b_{K_j}\Big] \\ \\ 
\ds -\lambda_i b_{\sum_{j=1}^I \lambda_j K_j}
 = -\lambda_i \mathrm{Cov}(\gamma_{K_i})b_{\sum_{j=1}^I \lambda_j K_j}.
\end{array}$$

\hfill \\
Since $\gamma_n(K_i)>0$, $\gamma_{K_i}$ cannot be supported in an affine hyperplane. Therefore, $\mathrm{Cov}(\gamma_{K_i})$ is positive definite and in particular invertible so that $$r_{K_i}(0)b_{K_i}=r_{\sum_{j=1}^I \lambda_j K_j}(0) \, b_{\sum_{j=1}^I \lambda_j K_j}.$$

\hfill \\
We now explain how equality in the uncentered Gaussian Brunn-Minkowski inequality forces convexity up to null sets. 
Let $C_i=\mathrm{conv}(K_i)$ for any $1 \le i \le I.$ Then, $$\sum_{i=1}^I \lambda_i C_i=\mathrm{conv}\Big(\sum_{i=1}^I \lambda_i K_i\Big)=\sum_{i=1}^I \lambda_i K_i:=K.$$ Consequently, replacing $K_i$ by any subset $K_i \subset L_i \subset C_i$ leaves its Minkowski sum equal to $K.$
Now, suppose that $\gamma_n(C_i\setminus K_i)>0$ for some $1 \le i \le I.$ We choose a bounded measurable set $A \subset C_i\setminus K_i$ of positive Gaussian measure and subsets $A_t \subset A$ such that $\gamma_n(A_t)=t \downarrow 0.$ Define $$L_i(t)=K_i \cup A_t, \quad \xi_t:=\frac1t \int_{A_t}x\, d\gamma_n(x),$$ for any $t>0$, uniformly bounded vectors. Hence, $$\gamma_n(L_i(t))=\gamma_n(K_i)+t, \quad b_{L_i(t)}=b_{K_i}+\frac{t}{\gamma_n(K_i)}(\xi_t-b_{K_i})+O(t^2),$$ with a uniform error. Indeed, $$\begin{array}{lll}
\ds b_{L_i(t)} &=& \ds \frac1{\gamma_n(L_i(t))}\int_{L_i(t)}x\, d\gamma_n(x)=\frac1{\gamma_n(K_i)+t}\big(\gamma_n(K_i)b_{K_i}+t \xi_t\big)=b_{K_i} \\ \\
\nm
&& \ds +\frac t{\gamma_n(K_i)}(\xi_t-b_{K_i})-\frac{t^2}{\gamma_n(K_i)(\gamma_n(K_i)+t)}(\xi_t-b_{K_i}).
\end{array}$$
By definition of $A$, there exists $R<+\infty$ such that $A \subset B_R(0)$ and so, $$|\xi_t| \le \frac1t \int_{A_t}|x|\, d\gamma_n(x) \le R.$$ This implies $$\Big| b_{L_i(t)}-b_{K_i}-\frac t{\gamma_n(K_i)}(\xi_t-b_{K_i})\Big| \le \frac{R+|b_{K_i}|}{\gamma_n(K_i)^2}t^2.$$
\hfill \\
If we expand the deficit $$\begin{array}{lll}
\ds D(t)&:=& \ds e^{|\sum_{j \ne i}\lambda_j b_{K_j}+\lambda_i b_{L_i(t)}|^2/2n}\gamma_n(K)^{1/n} \ds -\sum_{j \ne i} \lambda_j e^{|b_{K_j}|^2/2n}\gamma_n(K_j)^{1/n} \\ 
\nm 
&&\ds -\lambda_i e^{|b_{L_i(t)}|^2/2n}\gamma_n(L_i(t))^{1/n} \ge 0,
\end{array}$$ one obtains that $$D(t)=\frac{\lambda_i t}{n\gamma_n(K_i)}\Big[(r_K(0)b_K-r_{K_i}(0)b_{K_i})\cdot (\xi_t-b_{K_i})-r_{K_i}(0)\Big]+O(t^2)=-\frac{\lambda_i r_{K_i}(0)}{n\gamma_n(K_i)}t+O(t^2).$$
More precisely, this holds as $$\begin{array}{lll}
\ds e^{\Big|\sum_{j=1}^I \lambda_j b_{K_j}+\lambda_i(b_{L_i(t)}-b_{K_i}) \Big|^2 \Big/2n}\gamma(K)^{1/n} &=& \ds r_K(0)\exp\Big(\frac{\lambda_i}n b_K \cdot (b_{L_i(t)}-b_{K_i})+\frac{\lambda_i^2}{2n}|b_{L_i(t)}-b_{K_i}|^2\Big) \\ \\
\nm 
&=& \ds r_K(0)+\frac{\lambda_i r_K(0)}{n}b_K\cdot (b_{L_i(t)}-b_{K_i})+O(t^2).
\end{array}$$
In addition, $$\begin{array}{lll}
\ds e^{|b_{L_i(t)}|^2/2n}\gamma_n(L_i(t))^{1/n} &=& \ds r_{K_i}(0)\Big(1+\frac t{\gamma_n(K_i)}\Big)^{1/n}\exp\Big(\frac1n b_{K_i}\cdot(b_{L_i(t)}-b_{K_i})+\frac1{2n}|b_{L_i(t)}-b_{K_i}|^2 \Big) \\ \\
\nm 
&=& \ds r_{K_i}(0)+\frac{r_{K_i}(0)}{n}b_{K_i}\cdot (b_{L_i(t)}-b_{K_i})+\frac{r_{K_i}(0)}{n\gamma_n(K_i)}t+O(t^2).
\end{array}$$
\hfill \\
Consequently, the deficit would be negative for small $t$ contradicting \cref{Generalized_Gaussian_BM}. Therefore, $$\gamma_n(C_i\setminus K_i)=0 \quad \forall 1 \le i \le I.$$ 

\hfill \\
The analysis is now split into two cases. If $b_{\sum_{j=1}^I \lambda_j K_j}=0$, this implies that $b_{K_i}=0$ for any $1 \le i \le I$ or in other words, all $K_i$ are centered.
\hfill \\
From \cref{Complete_BBL_eq} applied to $p=1, \tilde p=+\infty,$ and $$\Gamma=\Sigma_i=I_n, \quad S_i=0, \quad f_i=\mathbbm{1}_{C_i} \quad  \forall 1 \le i \le I \text{ and } h=\mathbbm{1}_K,$$ we conclude that $$K_1=\cdots=K_I,$$ up to Gaussian null sets.  

\hfill \\
Else, $b_{\sum_{j=1}^I \lambda_j K_j} \ne 0.$ Then, $$b_{K_i}=\frac{r_{\sum_{j=1}^I \lambda_j K_j}(0)}{r_{K_i}(0)}b_{\sum_{j=1}^I \lambda_j K_j} \quad \forall 1 \le i \le I,$$ and so, $$1=\sum_{i=1}^I \lambda_i \frac{r_{\sum_{j=1}^I \lambda_j K_j}(0)}{r_{K_i}(0)}.$$ Moreover, by equality, $$\sum_{i=1}^I \lambda_i \frac{r_{K_i}(0)}{r_{\sum_{j=1}^I \lambda_j K_j}(0)}=1,$$ and hence, $$1=\sum_{i=1}^I \lambda_i \frac{r_{\sum_{j=1}^I \lambda_j K_j}(0)}{r_{K_i}(0)} \sum_{i=1}^I \lambda_i \frac{r_{K_i}(0)}{r_{\sum_{j=1}^I \lambda_j K_j}(0)} \ge \Big(\sum_{i=1}^I \lambda_i\Big)^2=1,$$ by Cauchy-Schwarz. Since there is equality, $$r_{K_i}(0)=r_{\sum_{j=1}^I \lambda_j K_j}(0) \quad \forall 1 \le i \le I,$$ and all barycenters also coincide. This implies that $$\gamma_n(K_1)=\cdots=\gamma_n(K_I)=\gamma_n\Big(\sum_{i=1}^I \lambda_i K_i\Big).$$ We may apply \cref{Gaussian_PL_EqCases} to the functions $f_i=\mathbbm{1}_{C_i}$ and $h=\mathbbm{1}_{\sum_{j=1}^I \lambda_j K_j}$ with the pointwise assumption following from $\sum_{i=1}^I \lambda_i C_i=\sum_{j=1}^I \lambda_j K_j$ and equality of integrals follows from $\gamma_n(C_i)=\gamma_n(K_i)=\gamma_n\big(\sum_{i=1}^I \lambda_i K_i\big).$
Hence, these indicator functions are proportional almost everywhere and as they are non-zero, they coincide almost everywhere. Therefore, $$\gamma_n\Big(K_i \Delta \Big(\sum_{j=1}^I \lambda_j K_j\Big) \Big)=0 \quad \forall 1 \le i \le I.$$
\end{proof}

Using first the same variational argument as for the equality cases of the uncentered Brunn-Minkowski inequality and then revisiting the inverse Brascamp-Lieb proof of \cite{Shohei_Hiroshi2}, we may obtain the following.

\begin{theorem} \label{Eq_Cases_GCI_Uncentered}
   Under the hypotheses of \cref{Uncentered_GCI}, equality holds if and only if $b=0$ and there exist pairwise orthogonal subspaces $E_1,\cdots,E_I$ such that $$E_i \subset \ker(\Sigma_0^{-1/2}\Sigma_i\Sigma_0^{-1/2}-I_n), \quad \gamma_n\Big(\Sigma_0^{-1/2}K_i \Delta (D_i+E_i^\perp) \Big)=0,$$ where $D_i \subset E_i$ is convex, has positive Gaussian measure in $E_i$ and satisfies $\mathrm{bar}_{\gamma_{E_i}}(D_i)=0$ for any $1 \le i \le I.$ 
\end{theorem}

\begin{proof}
    First, suppose \cref{Uncentered_GCI_ineq} holds with equality. Denote $$K_0=\bigcap_{i=1}^I K_i.$$
    As previously, we may apply \cite[Theorem 1.5]{Shohei_Hiroshi2} to the convex sets $C_N(K_i-u,b-u)=C_N(K_i,b)$ for any $1 \le i \le I$, so that $$e^{\frac12 \langle \Sigma_0^{-1}(b-u),b-u\rangle} \gamma_{\Sigma_0}(K_0-u) \ge \prod_{i=1}^I \Big[ e^{\frac12 \langle \Sigma_i^{-1}b,b\rangle} \gamma_{\Sigma_i}(K_i) \Big] \quad \forall u \in \R^n.$$  We may define $$F:\mathbb{R}^n \rightarrow [0;+\infty), \quad F(u)=e^{\frac12 \langle \Sigma_0^{-1}(b-u),b-u\rangle} \gamma_{\Sigma_0}(K_0-u).$$ Then, by assumption, $F(u) \ge F(0)$ for any $u \in \R^n$, so that by differentiating $\log F$, we get similarly to the previous proof that $$b_{K_0,\Sigma_0}=b.$$ Write $$\tilde K_i=\Sigma_0^{-1/2}K_i \quad \forall 0 \le i \le I, \quad \tilde b=\Sigma_0^{-1/2}b.$$ 
    The covariances are now $$\Gamma_0=I_n, \quad \Gamma_i=\Sigma_0^{-1/2}\Sigma_i\Sigma_0^{-1/2}, \quad S=\sum_{i=1}^I \Gamma_i^{-1}-I_n.$$ For $i=0,\cdots,I$, define $$f_i(x)=\frac{\mathbbm{1}_{\tilde K_i}(x+\tilde b) \gamma_{\Gamma_i}(x+\tilde b)}{\gamma_{\Sigma_i}(K_i)\gamma_{\Gamma_i}(x)}=\frac{e^{-\langle \Gamma_i^{-1}\tilde b,\tilde b\rangle/2}}{\gamma_{\Sigma_i}(K_i)}e^{-\langle \Gamma_i^{-1}\tilde b,x\rangle}\mathbbm{1}_{\tilde K_i-\tilde b}(x).$$ Then, $f_i$ are log-concave densities with $\int f_i\, d\gamma_{\Gamma_i}=1$ and barycenter $0.$ The assumed equality yields $$f_0(x)=e^{\langle S\tilde b,x\rangle}\prod_{i=1}^I f_i(x).$$  Denoting $f_i^{(k)}=\mathrm{Conv}_{\gamma_{\Gamma_i}}^k(f_i)$, we claim that \begin{equation} \label{Doubling_Trick_GCI} f_0^{(k)}(x) \ge e^{2^{k/2}\langle S\tilde b,x\rangle}\prod_{i=1}^I f_i^{(k)}(x).\end{equation} Indeed, by induction, fix $x$ and write $$H_{i,x}(y)=f_i^{(k)}\Big(\frac{x+y}{\sqrt 2}\Big)f_i^{(k)}\Big(\frac{x-y}{\sqrt 2}\Big),$$ which are even and log-concave functions. By \cite[Theorem 2.4 together with Lemma 4.1]{Shohei_Hiroshi2}, $$\int \prod_{i=1}^I H_{i,x}\, d\gamma_n \ge \prod_{i=1}^I \int H_{i,x}\, d\gamma_{\Gamma_i}.$$ Multiplying \cref{Doubling_Trick_GCI} at $(x+y)/\sqrt{2}$ and $(x-y)/\sqrt 2$ and applying the previous inverse Brascamp-Lieb inequality, gives the desired induction step.  
    \hfill \\
    Letting $A_i=\mathrm{Cov}(f_i\gamma_{\Gamma_i})$, we have $0<A_i \le \Gamma_i$ and the local central limit theorem gives $$f_i^{(k)}(x) \rightarrow_{k \rightarrow +\infty} G_i(x)=\frac{\gamma_{A_i}(x)}{\gamma_{\Gamma_i}(x)}.$$ If $S\tilde b \ne 0$, we may apply \cref{Doubling_Trick_GCI} at $x=S\tilde b$, so that the left-hand side would have a finite limit, whereas the right hand side would tend to infinity using that $G_i(S\tilde b)>0.$ Hence, $S\tilde b=0.$ 
    \hfill \\
    Taking the limit in \cref{Doubling_Trick_GCI} gives $G_0 \ge \prod_{i=1}^I G_i,$ and we may again apply the inverse Brascamp-Lieb inequality to obtain that $$1=\int G_0\, d\gamma_n \ge \int \prod_{i=1}^I G_i \, d\gamma_n \ge \prod_{i=1}^I \int G_i\, d\gamma_{\Gamma_i}=1.$$ By continuity, this implies that $$G_0=\prod_{i=1}^I G_i.$$ 
    \hfill \\
    Let us now investigate the matrix equality cases. Set $$B_i=A_i^{-1}-\Gamma_i^{-1} \ge 0 \quad \forall 0 \le i \le I.$$ Then, $B_0=\sum_{i=1}^I B_i$ and $$\det\Big(I_n+\sum_{i=1}^I B_i\Big)=\prod_{i=1}^I \det(I_n+B_i^{1/2}\Gamma_i B_i^{1/2}).$$ As $$\begin{array}{lll}
    \ds \det(I_n+A+B) &=& \ds \det(I_n+A)\det(I_n+B^{1/2}(I_n+A)^{-1}B^{1/2})
    \\ 
    \nm 
    &\le& \ds \det(I_n+A)\det(I_n+B),
    \end{array}$$ and equality is equivalent to $AB=0$, we have $$\det\Big(I_n+\sum_{i=1}^I B_i\Big) \le \prod_{i=1}^I \det(I_n+B_i) \le \prod_{i=1}^I \det(I_n+B_i^{1/2}\Gamma_iB_i^{1/2}).$$ In particular, equality implies that $$B_iB_j=0 \quad \forall i \ne j, \quad B_i^{1/2}(\Gamma_i-I_n)B_i^{1/2}=0 \quad \forall i,$$ or equivalently, $(\Gamma_i-I_n)B_i=0$ for any $1 \le i \le I.$ 
    
    \hfill \\
    We may define the pairwise orthogonal spaces $E_i=\mathrm{Ran}(B_i)$ for $1 \le i \le I$. In addition, $\Gamma_i|_{E_i}=I_{E_i}$ and by symmetry of $\Gamma_i$, $E_i^\perp$ is also an invariant subspace of $\Gamma_i$. As $A_i=(\Gamma_i^{-1}+B_i)^{-1}$, $$A_iv=\Gamma_i v \quad \forall v \in E_i^\perp.$$ 
    \hfill \\
    By \cref{Splitting_Lemma}, $$\frac{e^{-\langle \Gamma_i^{-1}\tilde b,\tilde b\rangle/2}}{\gamma_{\Sigma_i}(K_i)}e^{-\langle \Gamma_i^{-1}\tilde b,x\rangle}\mathbbm{1}_{\tilde K_i-\tilde b}(x)=f_i(x)=\psi_i(P_{E_i}x) \text{ for a.e. }x.$$ Considering the convex subset $D_i=\{\psi_i>0\} \subset E_i$, we get that $$\tilde K_i-\tilde b=D_i+E_i^\perp \text{ up to null sets}, \text{ and } \Gamma_i^{-1}\tilde b \in E_i \text{ for any } 1 \le i \le I,$$ as $x \mapsto \langle \Gamma_i^{-1}\tilde b,x\rangle$ is constant on $E_i^\perp$ and therefore identically $0.$
    Thus, $\tilde b \in E_i$ for any $1 \le i \le I$ but as $I \ge 2$, this implies that $\tilde b=0$, and so also $b=0.$ 
    In other words, $$\Sigma_0^{-1/2}K_i=\tilde K_i=D_i+E_i^\perp \quad \forall 1 \le i \le I,$$ with $D_i$ being centered with respect to $\gamma_{E_i}$, as the marginal of $\gamma_{\Gamma_i}$ on $E_i$ is the standard Gaussian on $E_i.$ 
    \hfill \\
    Conversely, if $b=0$, $$\gamma_{\Sigma_i}(K_i)=\gamma_{\Sigma_0^{-1/2}\Sigma_i\Sigma_0^{-1/2}}(\Sigma_0^{-1/2}K_i)=\gamma_{E_i}(D_i) \quad \forall 1 \le i \le I,$$ while $$\gamma_{\Sigma_0}\Big(\bigcap_{i=1}^I K_i\Big)=\gamma\Big(\bigcap_{i=1}^I (\Sigma_0^{-1/2}K_i)\Big)=\prod_{i=1}^I \gamma_{E_i}(D_i),$$ by independence, using that $E_1,\cdots,E_I$ are pairwise orthogonal subspaces. 
\end{proof}

\hfill \\

\section*{Acknowledgements}
The author is grateful to his supervisor Alessio Figalli, for his constant support and encouragement. He would also like to thank Hiroshi Tsuji for a useful comment about \cref{Counterexample_FailureLogCcv}. 

\hfill \\
\textbf{AI Usage:} All the proofs and ideas were developed by the author, except the proof of \cref{Counterexample_FailureLogCcv} proposed by a commercial AI tool, checked and rewritten by the author. 
\hfill \\
\printbibliography
 
\end{document}